\documentclass{siamart0516}
\usepackage[margin=1in]{geometry}
\usepackage{amsmath,comment}
\usepackage{amssymb}
\usepackage{amsfonts} 
\usepackage[normal]{subfigure}
\providecommand{\U}[1]{\protect\rule{.1in}{.1in}}

\newtheorem{exmp}{Example}[section]
\renewcommand\thesection{\arabic{section}}
\renewcommand\thesubsection{\thesection.\arabic{subsection}}
\renewcommand\thesubsubsection{\thesubsection.\arabic{subsubsection}}
\renewcommand\theequation{\thesection.\arabic{equation}}
\renewcommand{\eqref}[1]{(\ref{#1})}
\usepackage[numbers,sort&compress]{natbib}

\newcommand{\eps}{{\displaystyle \varepsilon}}
\newcommand{\bsub}{\begin{subequations}}
\newcommand{\esub}{\end{subequations}$\!$}

\newcommand{\littleoh}{ \mbox{{\scriptsize $\mathcal{O}$}}}
\newcommand{\bigoh}{\mathcal{O}}

\newcommand{\bx}{\mathbf{x}}

\newcommand{\ba}{\mathbf{a}}

\newcommand{\by}{\mathbf{y}}
\newcommand{\bz}{\mathbf{z}}
\newcommand{\bp}{\mathbf{p}}

\newcommand{\Q}{\mathcal{Q}}

\newcommand{\bxi}{\boldsymbol{\xi}}

\newcommand{\M}{{\mathcal{M}}}
\newcommand{\I}{{\mathcal{I}}}
\newcommand{\A}{{\mathcal{A}}}

\newcommand{\hn} {{\hat{\textbf{n}}}}
\begin{document}

\title{The effect of target orientation on the mean first passage time of a Brownian particle to a small elliptical absorber.}

\date{\today}
\author{Sanchita Chakraborty\thanks{Department of Applied and Computational Mathematics and Statistics, University of Notre Dame, Notre Dame, IN, 46656, USA. {\tt schakra2@nd.edu}}\and Theodore Kolokolnikov\thanks{Department of Mathematics, Dalhousie University, Halifax, Nova Scotia, Canada, B3H 3J5. {\tt tkolokol@gmail.com}}
\and Alan E. Lindsay\thanks{Department of Applied and Computational Mathematics and Statistics, University of Notre Dame, Notre Dame, IN, 46656, USA. {\tt a.lindsay@nd.edu}}
}

\maketitle

\begin{abstract}
We develop a high order asymptotic expansion for the mean first passage time (MFPT) of the capture of Brownian particles by a small elliptical trap in a bounded two dimensional region. This new result describes the effect that trap orientation plays on the capture rate and extends existing results that give information only on the role of trap position on the capture rate. Our results are validated against numerical simulations which confirm the accuracy of the asymptotic approximation. In the case of the unit disk domain, we identify a bifurcation such that the high order correction to the global MFPT (GMFPT) is minimized when the trap is orientated in the radial direction for traps centered at $0<r<r_c :=\sqrt{2-\sqrt{2}}$. When centered at position $r_c<r<1$, the GMFPT correction is minimized by orientating the trap in the angular direction. In the scenario of a general two-dimensional geometry, we identify the orientation that minimizes the GMFPT in terms of the regular part of the Neumann Green's function. This theory is demonstrated on several regular domains such as disks, ellipses and rectangles.
\end{abstract}


\section{Introduction}

We consider the problem of describing the mean first passage time (MFPT) of two dimensional Brownian motion in a bounded region to a small elliptical absorbing trap. The diffusive transport of molecules and individual agents from a source to a mobile or fixed target is a problem occurring in a variety of physical, biological and social systems \cite{REDNER2001,Schuss2014,Metzler2014}. Ecological examples include the time required for an animal to find a mate or shelter \cite{Venu2015,Stepien2020,FAUCHALD2003}. At the cellular scale, diffusion transports key cargoes within the cell \cite{Miles2020,Bressloff2024,GHUSINGA2017,GIBBONS2010,Grebenkov_2017}, including fibroblasts to initiate wound healing \cite{ambrosi2004cell}, antigens for detection by T-cell receptors \cite{Morgan2023,MorganLindsay2022}, and material to and from the nucleus \cite{WINDNER2019,LEECH2022}. The scope of the target search problem has also been expanded to include features such as stochastic switching of target states, resetting \cite{BressloffSchumm2022}, extreme statistics, and homogenization \cite{Lawley2023}. For an extensive review, we point the reader to the recent survey \cite{Grebenkov2024}.

The MFPT $u(\bx)$ describing the expected time to capture a diffusing particle initially at $\bx\in\Omega\setminus\Omega_{\eps}$ solves the Poisson equation \cite{Newby2013,Schuss2014}
\bsub\label{eqn:MFPTIntro}
\begin{align}
\label{eqn:MFPTIntro_a} D \Delta u+1 &= 0, \qquad \bx\in\Omega\setminus\Omega_{\eps};\\[4pt]
\label{eqn:MFPTIntro_b}  D\nabla u \cdot \hn &= 0, \qquad \bx \in \partial\Omega;\\[4pt]
\label{eqn:MFPTIntro_c} u &= 0, \qquad \bx \in \partial\Omega_{\eps}.
\end{align}
\esub
The boundary conditions \eqref{eqn:MFPTIntro_b} prescribe that the outer boundary $\partial\Omega$ is reflecting and equation \eqref{eqn:MFPTIntro_c} specifies that the trap $\Omega_{\eps}$ is absorbing. The aim of this paper is to construct a solution to \eqref{eqn:MFPTIntro} in the limit as $\eps\to0$ in the presence of an elliptical trap defined as 
\begin{equation}\label{eq:ellipse}
\Omega_{\eps} = \bxi + \eps e^{i\phi} \A, \qquad \A = \Big\{ (y_1,y_2) \in \mathbb{R}^2 \ | \ \frac{y_1^2}{a^2} + \frac{y_2^2}{b^2} < 1\Big\}. 
\end{equation}

\begin{figure}[htbp]
\centering
\includegraphics[width = 0.75\textwidth]{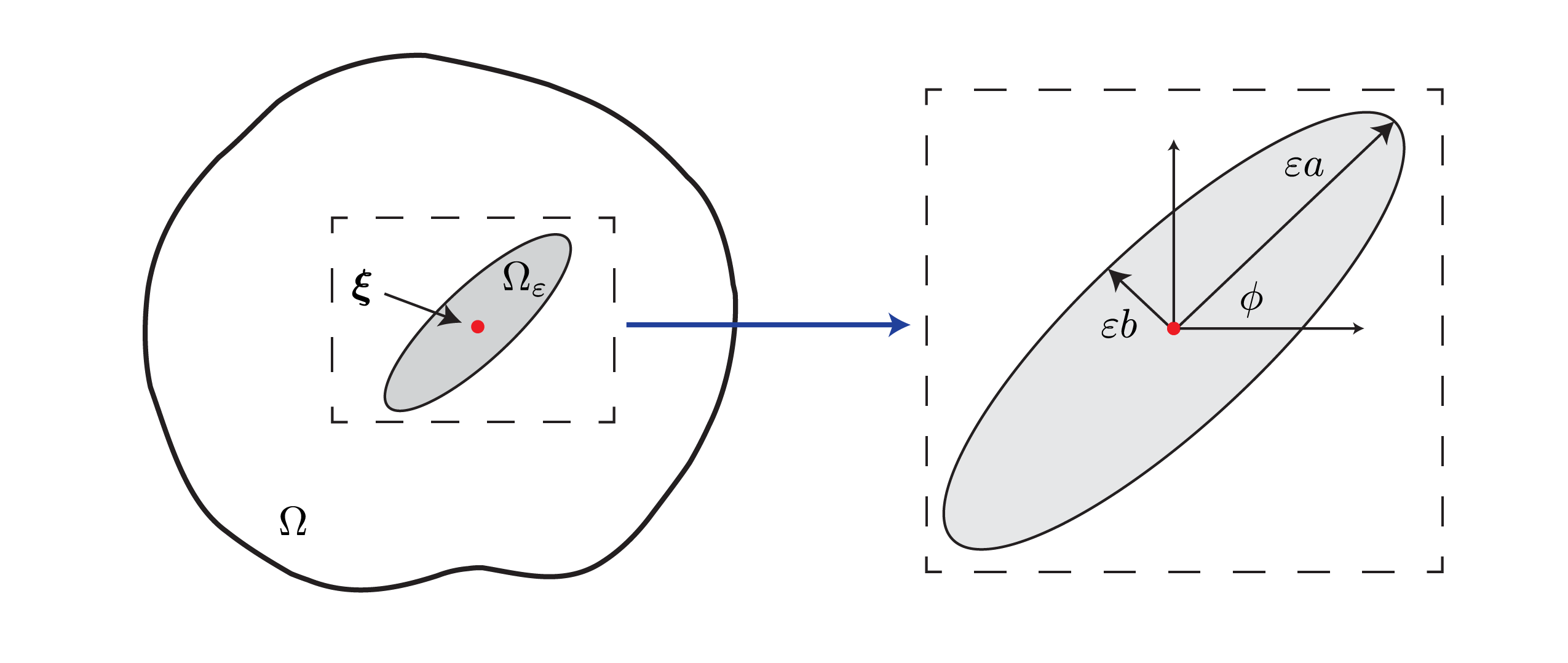}
\caption{Schematic of the configuration of the domain $\Omega$ with a single trap $\Omega_{\eps}$ as defined in \eqref{eq:trap:def}. The trap is centered at $\bxi\in\Omega$ and has semi-major and semi-minor axes $\eps a$ and $\eps b$ respectively. The semi-major axis of the trap is orientated at angle $\phi$ with respect to the horizontal axis. \label{fig:intro}}
\end{figure}

Here $\eps a$ and $\eps b$ are the semi-major and semi-minor axes respectively and $\phi$ is the angle of orientation with respect to the horizontal axis (see Fig.~\ref{fig:intro}). The term $e^{i\phi}$ corresponds to rotation by angle $\phi$ in the counter clockwise direction. An important quantity, called the global MFPT (GMFPT), describes the overall capture rate based on a uniform distribution of start locations, and is defined as
\begin{equation}\label{eq:trap:def}
\tau =\frac{1}{|\Omega\setminus\Omega_{\eps}|} \int_{\Omega\setminus\Omega_{\eps}} u(\bx)\, d\bx.
\end{equation}

Before outlining the rationale for this work and relationship to previous studies, we state our main result:

\textit{Principal Result: Consider equation \eqref{eqn:MFPTIntro} with a single elliptical trap centered at $\bxi\in\Omega$ with semi-major and semi-minor axes $\eps a$ and $\eps b$ respectively ($a>b$) and the semi-major axis having elevation $\phi$ from the horizontal. In the limit as $\eps\to0^{+}$, a two term expansion of the solution to
\eqref{eqn:MFPTIntro} and the GMFPT \eqref{eq:trap:def} is of form}
\bsub\label{eq:introExpansion}
\begin{equation}
\label{eq:introExpansion_a} u(\bx) = \frac{1}{D}\left[ u_0(\bx) + \eps^2 u_2(\bx) + \bigoh(\eps^4)\right]; \qquad  \tau = \frac{1}{D}\left[\tau_0 + \eps^2 \tau_2  + \bigoh(\eps^4)\right].
\end{equation}
\textit{The terms in the above expansions are given explicitly as }
\begin{align}
\label{eqn:intro_u0} u_0(\bx) &= -|\Omega|\Big[ G(\bx;\bxi) - R(\bxi;\bxi) \Big] + \frac{|\Omega|}{2\pi  \nu};\\[4pt]
\label{eqn:intro_u2} u_2(\bx) &=  |\Omega| \Big [\frac12 \mbox{Trace} \big( \Q \nabla^2_{\bxi} G(\bx;\bxi)\big) -2\pi\nabla_{\bxi}R(\bxi;\bxi) \cdot \M \nabla_{{\bxi}} G(\bx;\bxi) \Big] + \chi_2.
\end{align}
\textit{The logarithmic gauge function is $\nu(\eps)= -1/\log (\eps d_c)$ where $d_c$ is the \emph{logarithmic capacitance} which reflects the shape of the trap and is determined by \eqref{vc}. For an elliptical trap, $d_c = \frac{1}{2}(a+b)$. The constant $\chi_2$ is given by} 
\begin{equation}
    \chi_2 = -|\Omega|\Big( \mbox{Trace}\big(\Q \nabla_{{\bxi}}^2 R(\bxi;\bxi) \big) - 2\pi \nabla_{\bxi}R(\bxi;\bxi)\cdot \M \nabla_{\bxi}R(\bxi;\bxi) \Big) + \frac{a^2+b^2}{8}.
\end{equation}
\textit{Here $\Q$ and $\M$ are the quadrupole and moment polarization matrices and solve associated electrified disk problems \eqref{vc} and \eqref{vc1}, respectively. For the case of an elliptical trap, they are given explicitly by}
\begin{equation}
\Q = -\frac{a^2-b^2}{4} \begin{bmatrix} \cos2\phi& \phantom{-}\sin2\phi \\ \sin2\phi & -\cos2\phi \end{bmatrix}, \qquad \M = -\frac{(a+b)^2}{4}\I + \Q.
\end{equation}
\textit{The terms of the GMFPT \eqref{eq:introExpansion_a} are given by }
\begin{align}
\label{eqn:u0_tau0} \tau_0 & =  \frac{|\Omega|}{2\pi }\Big[ \frac{1}{\nu} + 2\pi R(\bxi;\bxi) \Big], \qquad \nu(\eps) = \frac{-1}{\log( \eps d_c)};\\
\label{eqn:u0_tau2} \tau_2 & =  \,\Big[ \frac{\pi ab }{|\Omega|}\tau_0 + \frac{a^2+b^2}{8} + \chi_2 \Big].
\end{align}
\esub

In the above result $G(\bx;\bxi)$ and $R(\bx;\bxi)$ are the Neumann Green's function and its regular part respectively, defined as the unique solution of
\bsub\label{neum_g}
\begin{gather}
\Delta G = \frac{1}{|\Omega|} - \delta(\bx-\bxi), \quad \bx\in\Omega; \qquad \nabla G \cdot \hn = 0,\quad \bx\in\partial\Omega;\\
 \int_{\Omega} G(\bx;\bxi) d\bx = 0; \qquad
G(\bx;\bxi) = -\frac{1}{2\pi}\log | \bx- \bxi| + R(\bx;\bxi).
\end{gather}
\esub

Before giving a detailed outline of the steps leading to the principal result, we review some recent work on related problems and motivations for this study. Over the past several decades there has been extensive research in the asymptotic analysis of the two-dimensional Poisson problem \eqref{eq:ellipse} in the presence of small inhomogeneities \cite{Iyaniwura2021_a,Iyaniwura2021_b,WK1993,KTW2005,bressloff2023asymptotic,Venu2015,NewbyIsaacson2013} which serves as a canonical problem in the trafficking and delivery of small signaling molecules and cargoes \cite{Newby2013}. More generally, there has been significant recent interest in the study of elliptic problems in punctured domains \cite{LLQ17,kropinski2011asymptotic,LHS2016,LWK2015,kolokolnikov2015recovering,Campbell1998}.

Elliptical traps are of particular interest in cellular signaling problems due to the frequent observation of non-circularity in the cell itself \cite{kaiyrbekov2024does} or other key organelles, such as the nucleus \cite{Manhart2025,WINDNER2019}. An oval or elliptical geometry accurately captures the aberrations to radial symmetry observed in these domains and hence it is natural to investigate how the capture rate of Brownian particles is modulated by non-circularity. This is one element in a broader mathematical effort to understand the contribution of geometry in cellular signaling \cite{nakamura2024gradient,ABNH2023,BERNOFF2023,kaiyrbekov2024does,Lindsay2025}

The leading order behavior of the GMFPT \eqref{eqn:u0_tau0} as $\eps\to0$ (see \cite{Venu2015}) captures the effects of trap size and position. However, there is no information on how the orientation of the trap influences the solution. A correction to the leading order behavior \eqref{eqn:intro_u0} was derived in \cite{Lindsay2015} that captures the effect of orientation,
\begin{equation}\label{eqn:intro_u1}
\tau = \frac{1}{D}\left[\tau_0 -  \eps|\Omega|\,\big( \, \textbf{d} \cdot \nabla_{\bxi}  R(\bxi;\bxi)\big) + \bigoh(\eps^2) \right].
\end{equation}
In the result \eqref{eqn:intro_u1}, the vector $\textbf{d}$, is related to the \emph{dipole moment} of the trap and is defined by associated problems that incorporate the shape and orientation of the trap. The contribution to the MFPT from the location of the trap in $\Omega$ is captured by the quantity $\nabla_{\bxi}  R(\bxi;\bxi)$ where the subscript reflects differentiation with respect to the source location.

For application of \eqref{eqn:intro_u1} to the case of an elliptical trap, which has two lines of symmetry, we find (see Appendix \ref{sec:Quadterm}) that the dipole term vanishes (\textbf{d}=\textbf{0}), thus \eqref{eqn:intro_u1} no longer describes the effect of orientation on the MFPT. Our refined result \eqref{eq:introExpansion} describes the higher order contribution to the MFPT due to the trap orientation. In particular, we can identify optimizing configurations by writing the GMFPT as
\begin{equation}\label{eqn:intro_tau}
\tau = \frac{\tau_0}{D} + \frac{\eps^2}{D} \left[\frac{\pi ab}{|\Omega|}\tau_0 + \frac{a^2+b^2}{4} - \pi |\Omega| \frac{(a+b)^2}{2} (R_{\xi_1}^2 + R_{\xi_2}^2) + |\Omega| \frac{a^2-b^2}{4} \bp \cdot \begin{bmatrix}\cos 2\phi \\ \sin 2\phi \end{bmatrix} \right].
\end{equation}
The vector $\bp$ is found to be
\begin{equation}\label{eqn:intro_optimal_p}
    \bp = \begin{bmatrix} 
R_{\xi_1\xi_1}-R_{\xi_2\xi_2} - 2\pi (R_{\xi_1}^2 -R_{\xi_2}^2) \\[4pt]
2R_{\xi_1\xi_2} - 4\pi R_{\xi_1} R_{\xi_2}
\end{bmatrix}\, ,
\end{equation}
which gives the direction along which the trap should be orientated to optimize the correction term $\tau_2$ of the GMFPT. We note from \eqref{eqn:u0_tau0}, that when the trap center $\bxi$ is placed at a critical point of $R(\bxi;\bxi)$, so that $\nabla_{\bxi}R(\bxi;\bxi) = [R_{\xi_1},R_{\xi_2}]^T=[0,0]^T$, the optimal orientation vector reduces to $\bp =[R_{\xi_1\xi_1}-R_{\xi_2\xi_2}, 2R_{\xi_1\xi_2}]^{T}$.

The outline of the paper is as follows. In Sec.~\ref{sec:trap_asy} we present a hierarchy of results, beginning with the solution of \eqref{eqn:MFPTIntro} in the reduced case of a circular trap located at the center of a disk (Sec.~\ref{sec:trap_asy_disk}). Following this, we derive the solution in the presence of an elliptical trap placed at the center of a disk (Sec.~\ref{sec:trap_asy_ellipse}). Finally, we present the corresponding result for the case of a general domain with an ellipse of arbitrary orientation (Sec.~\ref{sec:trap_asy_general}) which yields the principal result \eqref{eq:introExpansion}.

In Sec.~\ref{sec:results} we first validate the asymptotic result on the unit disk domain where the regular part $R(\bx;\bxi)$ is known in closed form. In this unit disk case, we identify a bifurcation where for $|\bxi|>r_c:=\sqrt{2-\sqrt{2}}$, the GMFPT correction term $\tau_2$ is minimized when the semi-major axis is orientated in the angular direction. Conversely, for $|\bxi|<r_c$, the GMFPT correction $\tau_2$ is minimized when the semi-major axis is orientated in the radial direction. For certain regular domains, such as rectangles and ellipses, highly accurate series solutions for \eqref{neum_g} are available, which allows us to determine $\bp$. For these cases we reveal similar bifurcations of the optimizing orientation depending on the centering point of the ellipse $\bxi$ and proximity to $\partial\Omega$. Finally in Sec.~\ref{sec:discusson}, we discuss avenues for future research arising from this study.

\section{Asymptotic analysis of the mean first passage time to a single elliptical trap}\label{sec:trap_asy}

In this section we perform the main asymptotic analysis on the MFPT problem \eqref{eqn:MFPTIntro}. To guide the rationale for the higher order expansions, it is useful to first analyze two exactly solvable cases for a circular trap located at the center of a disk and an elliptical trap located at the center of a disk.

\subsection{Unit disk with a circular trap at the origin}\label{sec:trap_asy_disk}
For a single circular trap at the origin, the MFPT \eqref{eqn:MFPTIntro} reduces to the ODE
\bsub
\begin{equation}\label{eqn:radialMPFT}
u_{rr} + \frac{1}{r} u_r = -\frac{1}{D}, \qquad \eps<r<1; \qquad u(\eps) = u'(1) = 0,
\end{equation}
where $r = |\bx|$. The exact solution of \eqref{eqn:radialMPFT} is
\begin{equation}\label{eq:exact_radial}
u(r) = \frac{1}{2D} \left[ - \frac{r^2}{2} + \log \frac{r}{\eps} + \frac{\eps^2}{2} \right].
\end{equation}
The corresponding GMFPT
\begin{align}\label{eq:globalMFPT_radial}
\nonumber \tau =& \frac{1}{|\Omega\setminus\Omega_{\eps}|} \int_{\Omega\setminus\Omega_{\eps}} u \, d\bx = \frac{2\pi}{\pi(1-\eps^2)} \int_{r=\eps}^1 u(r) r dr = \frac{1}{8D(1-\eps^2)} \Big[ -3 + 4\eps^2 - \eps^4 - 4\log \eps \Big]\\[4pt]
& = \frac{1}{8D} \Big[ -3 - 4 \log \eps + \eps^2 (1 - 4 \log\eps) + \bigoh(\eps^4) \Big].
\end{align}
\esub
The equation \eqref{eq:globalMFPT_radial} will be a useful case to validate solutions of \eqref{eqn:MFPTIntro} for more general configurations.

\subsection{Unit disk with an elliptical trap at the origin}\label{sec:trap_asy_ellipse}

We now solve for the MFPT in the scenario of an elliptical trap at the origin, orientated along the horizontal axis ($\phi=0$) with semi-major and semi-minor axes $\eps a$ and $\eps b$ respectively.

\underline{Outer Expansion} We expand as 
\begin{equation}\label{eq:expansion_ex2}
u(r) = \frac{1}{D} \left[ u_0(r) + \eps^2 u_2(r,\theta) + \cdots\right]
\end{equation}
where $\Delta u_0 +1 = 0$ with $u_0'(1) = 0$ and $\Delta u_2 = 0$ with $u_2'(1) = 0$. The general solutions for $u_0$ and $u_2$ in polar coordinates $\bx = re^{i\theta}$ are
\bsub\label{eq:CenterLeading}
\begin{align}
\label{eq:CenterLeading_a} u_0 &= \frac{1}{2}\log r - \frac{r^2}{4} + A_1\\[4pt]
\label{eq:CenterLeading_b} u_2 &= B_2 \cos 2\theta ( r^2 + r^{-2} ) + B_1,
\end{align}
\esub
for constants $A_1,B_1, B_2$ to be determined. We remark that in \eqref{eq:CenterLeading_b}, the general solution can accommodate a  term of form $\sin2\theta\times\{r^2,r^{-2}\}$ in the case that the trap orientation moves off the perpendicular axis. Moving to a coordinate $\bx = \eps \by$, we find that
\[
u \sim \frac{1}{2} \log |\by| + \frac{1}{2} \log\eps + A_1 + B_2 \frac{\cos 2\theta}{|\by|^2} + \eps^2 \Big[ B_1 -\frac{y^2}{4}  \Big] + \bigoh(\eps^4).
\]
\underline{Inner Expansion:} In the region $\bx = \eps\by$, the solution $u(\bx) = U(\by)$ is expanded as
\[
U(\by) = \frac{1}{D} \left[U_0(\by) + \eps^2 U_2(\by) + \cdots \right].
\]
\underline{$\bigoh(\eps^0):$} The leading order problem satisfies
\bsub\label{eq:eqnU0}
\begin{gather}
\label{eq:eqnU0_a} \Delta_{\by} U_0 = 0, \quad \mbox{in} \quad \mathbb{R}^2\setminus \A; \qquad U_0 = 0 \quad \mbox{on} \quad \partial\A;\\[4pt]
\label{eq:eqnU0_b} U_0 = \frac{1}{2} \log|\by| + \frac{1}{2} \log\eps + A_1 + B_2 \frac{\cos 2\theta}{|\by|^2} + \cdots, \quad |\by|\to\infty.
\end{gather}
where $\A$ is the rescaled ellipse of semi-major and semi-minor axes $a$ and $b$ respectively with orientation $\phi=0$ with respect to the origin. Using the solution $v_{0c}(\by)$ of the electrified disk problem derived in appendix \ref{sec:Quadterm}, we calculate that 
\begin{equation}\label{eqnU0farfield}
U_0 = \frac{1}{2}v_{0c}(\by) \sim \frac{1}{2} \Big[\log|\by| - \log \alpha - \frac{a^2-b^2}{4} \frac{\cos 2\theta}{|\by|^2} + \cdots \Big], \qquad |\by| \to \infty,
\end{equation}
\esub
where $\alpha = (a+b)/2$. Matching \eqref{eq:eqnU0_b} with \eqref{eqnU0farfield} yields that
\begin{equation}
A_1 = \frac{1}{2\nu}, \qquad B_2 = -\frac{a^2-b^2}{8}, \qquad  \nu = \frac{-1}{\log \eps\alpha}.
\end{equation}
\underline{$\bigoh(\eps^2):$} Proceeding to the next order, we have that
\bsub\label{eq:eqnU2}
\begin{gather}
\label{eq:eqnU2_a}  \Delta_{\by} U_2 = -1, \quad \mbox{in} \quad \mathbb{R}^2\setminus \A; \qquad U_2 = 0 \quad \mbox{on} \quad \partial\A;\\[4pt]
\label{eq:eqnU2_b} U_2 = B_1 - \frac{|\by|^2}{4} + \littleoh(1), \quad |\by|\to\infty.
\end{gather}
\esub
We solve this problem in appendix \ref{app:eqnv2} by decomposing the solution as $U_2 = -\frac{1}{4}|\by|^2 + U_{2h}$ where $U_{2h}(\by)$ solves a homogeneous problem. We find (see \eqref{app:2ndorder2}) that the large argument behavior of \eqref{eq:eqnU2} is
\begin{equation}
U_2 \sim - \frac{|\by|^2}{4} + \frac{a^2+b^2}{8} + \bigoh(|\by|^{-2}), \qquad |\by|\to\infty.
\end{equation}
Hence from comparison to \eqref{eq:eqnU2_b}, we determine that
\[
B_1 = \frac{a^2+b^2}{8}.
\]
We remark that the term $\mathcal{B}_{11}$ in \eqref{app:2ndorder2} is determined by the Hessian of the leading outer solution $U_0$ which vanishes in this case due to the trap being centered at the origin. This completes the derivation of the expansion \eqref{eq:expansion_ex2} and yields the following expression for the MFPT
\begin{equation}
u(r,\theta) =  \frac{1}{2D}\left[\log r - \frac{r^2}{2} + \frac{1}{\nu} +   \frac{\eps^2}{4} \Big( (a^2+b^2) -(a^2-b^2)\big(r^2+r^{-2} \big)\cos 2\theta\Big) \right] + \cdots
\end{equation}
as $\eps\to0$. We remark that this expression reduces to the exact solution \eqref{eq:exact_radial} in the scenario $a=b=1$. Calculating the GMFPT, we must identify local and global contributions by introducing an intermediate scale $\eps \ll \delta \ll1$.
\[
\tau = \frac{1}{|\Omega\setminus\Omega_{\eps}|} \int_{\Omega\setminus\Omega_{\eps}} u \, d\bx = \frac{1}{\pi(1- \eps^2ab)} \Big[ \underbrace{\int_{|\bx| <\delta} u d\bx}_{I_1} + \underbrace{\int_{\delta<|\bx|<1} u d\bx}_{I_2}\Big].
\]
We proceed to calculate both terms $I_1$ and $I_2$ to arrive at a final expression independent of $\delta$.

\underline{Calculation of $I_1$:} Rescaling with $\bx = \eps\by$, we have that
\begin{equation}
I_1 = \eps^2 \int_{\substack{\by\in\mathbb{R}^2\setminus\A\\ |\by|<\delta/\eps}} U\, d\by = \frac{\eps^2}{D} \int_{|\bz|=1}^{|\bz| = \delta/\eps} (\underbrace{U_0}_{I_{10}} + \eps^2 \underbrace{U_1}_{I_{12}})|J|d\bz.
\end{equation}
In the above calculation of the integral $I_1$, the region exterior to the ellipse has been mapped to the exterior of the disk through the transformation $\by\to\alpha \bz + \beta/\bz$ for $\bz = re^{i\theta}$ and Jacobian $|J|$ given by
\[
|J| = \begin{vmatrix} \alpha - \frac{\beta}{r^2}\cos 2\theta & - \frac{\beta}{r^2} \sin 2\theta \\[4pt] \frac{\beta}{r^2} \sin 2\theta & \alpha - \frac{\beta}{r^2}\cos 2\theta \end{vmatrix} = \alpha^2- \frac{2\alpha\beta}{r^2}\cos 2\theta + \frac{\beta^2}{r^4}.
\]
Using the fact that $U_0(\bz) = \frac{1}{2} \log |\bz|$, we have that
\begin{align}
\nonumber I_{10} &\approx \frac{\eps^2}{2} \int_{\theta=0}^{2\pi}\int_{r=1}^{r= \frac{\delta}{\eps\alpha}} \log r \Big(\alpha^2 - \frac{2\alpha\beta}{r^2} \cos 2\theta + \frac{\beta^2}{r^4} \Big) r dr d\theta = \eps^2 \pi \int_{r=1}^{r= \frac{\delta}{\eps\alpha}} \log r \Big(\alpha^2 r + \frac{\beta^2}{r^3} \Big) dr \\
& \approx  \frac{\pi}{4} \left[ \eps^2(\alpha^2+\beta^2)  + 2 \delta^2 \log \frac{\delta}{\eps \alpha} - \delta^2  \right].
\end{align}
The contribution from $I_{12}$ is $\bigoh(\eps^4)$. We now calculate the contribution from the outer solution.

\underline{Calculation of $I_2$:} 
\begin{align}
\nonumber I_2 &= \frac{1}{D}\int_{|\bx|>\delta} (u_0 + \eps^2 u_2) d\bx \\
\nonumber &=\frac{1}{2D} \int_{\theta=0}^{2\pi}\int_{r=\delta}^{1}\left[\log r - \frac{r^2}{2} + \frac{1}{\nu} +   \frac{\eps^2}{4} \Big( (a^2+b^2) -(a^2-b^2)(r^2 +r^{-2})\cos 2\theta \Big) \right] r dr d\theta\\
\nonumber {} & = \frac{\pi}{D}\int_{r=\delta}^{1}\left[ \log r - \frac{r^2}{2} + \frac{1}{\nu} +   \frac{\eps^2}{4} \Big(a^2+b^2\Big) \right] rdr \\
{}& \approx\frac{\pi}{4D}\left[ -\frac{3}{2} + \frac{2}{\nu} + \frac{\eps^2}{2}\Big(a^2+b^2\Big) -2 \delta^2 \log \delta + \delta^2 - \frac{2\delta^2}{\nu}  \right]
\end{align}
Combining these two terms, we have that 
\begin{align}
\nonumber \tau & = \frac{1}{\pi D(1-\eps^2 ab)} \frac{\pi}{8} \left[ -3 + \frac{4}{\nu} + \eps^2(a^2+b^2) + 2\eps^2 (\alpha^2+\beta^2)\right]\\[4pt]
\label{eqn:globaltau} & = \frac{1}{8D} \left[ -3 + \frac{4}{\nu} + \eps^2 \Big( 2(a^2+ b^2) - 3 ab + \frac{4ab}{\nu} \Big)   \right] + \bigoh(\eps^4).
\end{align}
As required, this expression is independent of $\delta$ and substitution of $a=b=1$ reduces \eqref{eqn:globaltau} to \eqref{eq:globalMFPT_radial}.

\subsection{General case for a single trap}\label{sec:trap_asy_general}
We will now determine the solution to the MFPT problem for a single elliptical domain, centered at $\bx=\bxi$ with semi-major and semi-minor axes $\eps a,\eps b$, respectively and orientation $\phi$ with respect to the horizontal axis. The explicit form of the trap is given in \eqref{eq:ellipse}.

In an outer region away form the elliptical trap, we expand the solution as  
\[
u(\bx) = \frac{1}{D} \left[u_0(\bx) + \eps u_1(\bx) + \eps^2 u_2(\bx) + \littleoh(\eps^2)\right].
\]
The outer problems $u_j$ for $j= 0,1,2,\ldots$ satisfy  
\bsub\label{eqn:outer_main}
\begin{align}
\Delta u_j + \delta_{1j} &= 0, \quad \bx\in\Omega\setminus \{\bxi\};\\[4pt]
\nabla u_j \cdot \hn &= 0, \quad \bx\in\partial\Omega.
\end{align}
\esub
The local behavior as $\bx\to\bxi$ is now established for each problem \eqref{eqn:outer_main} through boundary layer analysis. In the vicinity of the trap, the solution is expanded in variables
\begin{equation}\label{eqn:expInner}
U = \frac{1}{D} \left[ U_0(\by) + \eps U_1 (\by) + \eps^2 U_2(\by) + \littleoh(\eps^2)\right], \qquad \by = e^{-i\phi} \, \frac{\bx-\bxi}{\eps}.
\end{equation}
Collecting terms at relevant orders gives a sequence of problems to be solved.

\underline{Inner Region $\bigoh(\eps^0)$}: The leading order problem for $U_0$ satisfies
\bsub\label{eqn:OuterGen}
\begin{align}
\Delta U_0 &= 0, \quad \by\in\mathbb{R}^2\setminus \A;\\[5pt]
 U_0 &= 0, \quad \by \in \partial\A.
 \end{align}
\esub
In terms of the solution $v_{0c}$ of the electrified disk problem obtained in Appendix \ref{sec:Quadterm}, we have that
\begin{equation}\label{eqn:outerS}
U_0(\by) = S\nu v_{0c}(\by), \qquad \nu = \frac{-1}{\log\eps\alpha}.
\end{equation}
Here $S$ is a constant to be determined in the matching process. The far field of equation \eqref{eqn:outerS} supplies the appropriate local behavior for the outer solution. In \eqref{vc}, we establish that for an elliptical trap aligned in the horizontal direction ($\phi=0$), the far field behavior is
\begin{equation}
v_{0c}(\by) = \log |\by| - \log\alpha + \frac{\by^T \tilde{\Q} \by }{|\by|^4}, \qquad  \tilde{\Q} = - \alpha \beta\begin{bmatrix} 1 &\phantom{-}0 \\0& -1 \end{bmatrix} = - \frac{a^2-b^2}{4}\begin{bmatrix} 1 &\phantom{-}0 \\0& -1 \end{bmatrix}.
\end{equation}
Hence, in terms of the outer coordinate $\by = \eps^{-1} e^{-i\phi} (\bx-\bxi)$ incorporating rotation by $\phi$ with respect to the horizontal, and incorporating the far field behavior for $|\by|\to\infty$, we generate the local behavior as $\bx\to\bxi$
\begin{equation}\label{eqn:farfieldU0}
u\sim S\nu \left[ \log|\bx-\bxi|  + \frac{1}{\nu} + \eps^2\frac{(\bx-\bxi)^T \Q(\bx-\bxi) }{|\bx-\bxi|^4}\right] + \cdots.
\end{equation}
The matrix $\Q=e^{i\phi} \tilde{\Q}e^{-i\phi}$ is calculated as
\begin{equation}\label{eqn:QuadrapoleMatrix}
\Q = -\frac{a^2-b^2}{4}\begin{bmatrix} \cos\phi & -\sin\phi\\ \sin\phi & \phantom{-}\cos\phi  \end{bmatrix}\begin{bmatrix} 1& \phantom{-}0 \\ 0 & -1 \end{bmatrix}\begin{bmatrix} \phantom{-}\cos\phi & \sin\phi\\ -\sin\phi & \cos\phi  \end{bmatrix} = -\frac{a^2-b^2}{4} \begin{bmatrix} \cos2\phi & \phantom{-}\sin2\phi\\ \sin2\phi & -\cos2\phi \end{bmatrix}.
\end{equation}
This behavior is used to furnish terms in the outer expansion.

\underline{Outer Region $\bigoh(\eps^0)$}: The problem at leading order is
\bsub\label{eqn:sing_lead}
\begin{align}
\label{eqn:sing_lead_a} \Delta u_0 + 1 &= 0, \qquad \bx \in \Omega\setminus\{\bxi\};\\[5pt]
\nabla  u_0 \cdot \hn &= 0, \qquad \bx\in\partial\Omega; \\[5pt]
\label{eqn:sing_lead_b} u_0 &\sim S\nu \log|\bx-\bxi| + S  + \cdots, \qquad \bx\to\bxi.
\end{align}
\esub
In terms of the Neumann's Green's function \eqref{neum_g}, we have that
\begin{equation}\label{eqn:sing_sol}
u_0 = -2\pi S \nu G(\bx;\bxi) + \tau_0.
\end{equation}
where $\tau_0$ is a constant. The expansion of $u_0$ as $\bx\to\bxi$ gives the local behavior
\begin{align}\label{eqn:LOcorr}
\nonumber u_0 &= S\nu[\log|\bx-\bxi| - 2\pi R(\bx;\bxi)] + \tau_0\\[4pt]
 & \sim S\nu \log |\bx-\bxi| + \tau_0 - 2\pi S \nu \Big[  R(\bxi;\bxi) +\ba \cdot (\bx -\bxi ) + (\bx-\bxi)^T \nabla_{\bx}^2 R\mid_{\bx=\bxi} (\bx-\bxi) + \bigoh(|\bx-\bxi|^3) \Big],
\end{align}
where the coefficient terms for $\bx=(x_1,x_2)$ are given by
\begin{equation}
 \ba : = \begin{bmatrix} \partial_{x_1} R(\bx;\bxi)\\ \partial_{x_2} R(\bx;\bxi) \end{bmatrix}_{\bx = \bxi}, \qquad \nabla_{\bx}^2 R\mid_{\bx=\bxi} =  \frac{1}{2}\begin{bmatrix} \partial_{x_1x_1}R(\bx;\bxi) & \partial_{x_1x_2}R(\bx;\bxi)  \\ \partial_{x_1x_2}R(\bx;\bxi) & \partial_{x_2x_2}R(\bx;\bxi) \end{bmatrix}_{\bx=\bxi} = \frac12 \begin{bmatrix} R_{11} & R_{12}\\ R_{12} & R_{22} \end{bmatrix}.
\end{equation}

A system of two equations for unknowns $(S,\tau_0)$ is found by both matching \eqref{eqn:sing_sol} with the local behavior, and integrating \eqref{eqn:sing_lead_a}. This yields that
\begin{equation}\label{eqn:leading_order_sol}
S = \frac{|\Omega| }{2\pi\nu}, \qquad \tau_0 = \frac{|\Omega|}{2\pi\nu} \Big[1 + 2\pi\nu R(\bxi;\bxi) \Big].
\end{equation}

This result was established in \cite{KTW2005}. We now determine the correction to the inner expansion. In the local variable $\by = \eps^{-1}e^{-i\phi}(\bx-\bxi)$, equation \eqref{eqn:LOcorr} yields the far field behavior
\begin{equation}\label{eqn:outerfull}
U(\by) \sim S\nu \left[ \log|\by| - \log \alpha - 2\pi \eps\, \ba \cdot e^{i\phi} \by -2\pi \eps^2 \, \by^{T} e^{-i\phi} \nabla_{\bx}^2 R\mid_{\bx=\bxi} e^{i\phi} \by + \cdots \right], \qquad |\by| \to\infty.
\end{equation}
This reveals the leading order matching behavior for the higher order inner corrections. Specifically, we have that
\bsub
\begin{align}
U_1(\by) &\sim -2\pi S \nu \ba \cdot e^{i\phi} \by + \cdots\\[4pt]
U_2(\by) &\sim  -2\pi S \nu \by^{T} e^{-i\phi} \nabla_{\bx}^2 R\mid_{\bx=\bxi} e^{i\phi} \by + \cdots
\end{align}
\esub
as $|\by|\to\infty$. This behavior is now matched to corresponding inner problems.

\underline{Inner Region $\bigoh(\eps^1)$}: At this order we must solve the exterior problem 
\bsub\label{eq:U1_main}
\begin{gather}
\label{eq:U1_main_a} \Delta_\by {U}_{1} = 0\,, \quad  \by\in\mathbb{R}^2\setminus \A; \qquad  U_{1}  = 0, \quad \by \in \partial{\A};
\\[5pt]
 \label{eq:U1_main_b} U_{1} = -2\pi S \nu \ba \cdot e^{i\phi} \by + \cdots, \qquad |\by|\to\infty;
\end{gather}
\esub
In appendix \ref{app:orderone}, we introduce and solve the vector valued electrified disk problem $\textbf{v}_{1c}$
\bsub\label{eq:vc1_main}
\begin{gather}
\label{eq:vc1_main_a} \Delta_\by \textbf{v}_{1c} = 0\,, \quad  \by\in\mathbb{R}^2\setminus \A; \qquad  \textbf{v}_{1c}  = 0, \quad \by \in \partial\A;
\\[5pt]
 \label{eq:vc1_main_b} \textbf{v}_{1c} = \by + \frac{ \tilde{\M} \by}{|\by|^2} + \cdots \, \quad   |\by|\to\infty; \qquad \tilde{\M} = - \alpha\begin{pmatrix} a & 0\\ 0& b\end{pmatrix}.
\end{gather}
\esub
In terms of this solution, we write that 
\begin{equation}\label{eqn:U1}
U_1 = -2\pi S \nu \, \ba \cdot e^{i\phi} \textbf{v}_{1c}(\by).
\end{equation}
Applying the far field behavior \eqref{eq:vc1_main_b} as $|\by|\to\infty$, while returning to outer coordinates with variable $\by = \eps^{-1} e^{-i\phi} (\bx-\bxi)$, we generate the local behavior
\bsub\label{eqn:farfieldU1}
\begin{equation}
U_1 \sim -\frac{2\pi S \nu}{\eps} \, \ba \cdot \left[ (\bx-\bxi)  + \eps^2 \frac{\M (\bx-\bxi) }{|\bx-\bxi|^2} \right] + \cdots, \qquad \mbox{as} \quad \bx\to\bxi,
\end{equation}
In the case of an elliptical trap, we calculate explicitly that $\M = e^{i\phi} \tilde{\M} e^{-i\phi} $ has the form
\begin{align}
\nonumber \M &= -\alpha\begin{bmatrix} \cos\phi & -\sin\phi\\ \sin\phi & \phantom{-}\cos\phi  \end{bmatrix}\begin{bmatrix} a& 0 \\ 0 & b \end{bmatrix}\begin{bmatrix} \phantom{-}\cos\phi & \sin\phi\\ -\sin\phi & \cos\phi  \end{bmatrix} \\[4pt]
\nonumber &=-\alpha \begin{bmatrix} a\cos^2\phi + b\sin^2\phi & (a-b)\sin\phi\cos\phi \\ (a-b)\sin\phi\cos\phi & b\cos^2\phi + a\sin^2\phi  \end{bmatrix}\\[4pt]
\nonumber & = -\alpha^2\begin{bmatrix} 1& 0 \\ 0 & 1 \end{bmatrix}  -\alpha\beta\begin{bmatrix} \cos2\phi& \phantom{-}\sin2\phi \\ \sin2\phi & -\cos2\phi \end{bmatrix} \\[4pt]
& = -\alpha^2 \mathcal{I} + \Q; \label{eqn:M}
\end{align}
\esub
where $\Q$ is the quadrupole matrix derived in \eqref{eqn:QuadrapoleMatrix}.

\underline{Inner Region $\bigoh(\eps^2)$}: At $\bigoh(\eps^2)$ we introduce the scalar valued problem $U_2(\by)$ where
\bsub\label{eqnmain:V2}
\begin{gather}
\label{eqnmain:V2_a} \Delta_{\by} U_2  = -1, \quad  \by\in\mathbb{R}^2\setminus \A; \qquad U_2 = 0, \quad \by \in \partial \A; \\[5pt]
\label{eqnmain:V2_b} U_2 = -2\pi S\nu \, \by^T e^{-i\phi} \nabla_{\bx}^2 R\mid_{\bx=\bxi} e^{i\phi} \by + \cdots \, \quad   |\by|\to\infty;
\end{gather}
\esub
To further decompose the far field behavior \eqref{eqnmain:V2_b}, we write
\begin{align*}
     \mathcal{H}&=e^{-i\phi} \nabla_{\bx}^2 R\mid_{\bx=\bxi} e^{i\phi}\\[4pt]     
     &= \frac{1}{2}\begin{bmatrix} \phantom{-}\cos\phi & \sin\phi\\ -\sin\phi & \cos\phi  \end{bmatrix}\begin{bmatrix} R_{11}& R_{12} \\ R_{12} & R_{22} \end{bmatrix}\begin{bmatrix} \cos\phi & -\sin\phi\\ \sin\phi & \phantom{-}\cos\phi  \end{bmatrix}\\[4pt]
    & = \frac{R_{11} + R_{22}}{4} \begin{bmatrix} 1 & 0\\ 0 & 1 \end{bmatrix} + \frac{1}{4}\begin{bmatrix}
        (R_{11}-R_{22})\cos2\phi + 2 R_{12}\sin2\phi & \phantom{-}2R_{12}\cos2\phi - (R_{11}-R_{22}) \sin2\phi\\
        2R_{12}\cos2\phi - (R_{11}-R_{22})\sin2\phi & -(R_{11}-R_{22})\cos2\phi - 2 R_{12}\sin2\phi
    \end{bmatrix}\\[4pt]
    &=\frac{R_{11} + R_{22}}{4} \begin{bmatrix} 1 & 0\\ 0 & 1 \end{bmatrix} + \begin{bmatrix}
        \mathcal{B}_{11} & \phantom{-}\mathcal{B}_{12}\\
        \mathcal{B}_{12} & - \mathcal{B}_{11}        
    \end{bmatrix}.
\end{align*}
Here $\mathcal{B}$ is the matrix satisfying $\mbox{Trace}(\mathcal{B}) = 0$ with components
\bsub\label{eqn:MatrixB}
\begin{align}
\label{eqn:MatrixB_a}  \mathcal{B}_{11} &= \frac{1}{4} \Big[ (R_{11}-R_{22})\cos2\phi + 2 R_{12}\sin2\phi \Big],\\[4pt]
\label{eqn:MatrixB_b}  \mathcal{B}_{12} &= \frac{1}{4}\Big[ 2R_{12}\cos2\phi - (R_{11}-R_{22})\sin2\phi\Big].
\end{align}
\esub
The far field behavior is now in the form
\begin{equation}\label{eq:u2Far}
\by^T \mathcal{H} \by = \left[ (R_{11} +  R_{22})\frac{y_1^2 + y_2^2}{4} + \mathcal{B}_{11} (y_1^2 - y_2^2) + 2 \mathcal{B}_{12} \, y_1 y_2  \right].
\end{equation}
We remark that $\Delta u_0 = \Delta( S\nu \log |\bx-\bxi| - 2\pi S\nu R(\bx;\bxi) +\tau_0) = -1$ which implies that
\begin{equation}
\Delta R = R_{11} + R_{22} =  \frac{1}{|\Omega|}.
\end{equation}
After applying these reductions, together with $2\pi S\nu = |\Omega| $ from \eqref{eqn:leading_order_sol}, we can restate equation \eqref{eqnmain:V2} as 
\bsub\label{eqnmain:V2_new}
\begin{gather}
\label{eqnmain:V2_newa}\Delta_{\by} U_2 = -1, \quad  \by\in\mathbb{R}^2\setminus \A; \qquad U_2 = 0, \quad \by \in \partial \A;\\[5pt]
\label{eqnmain:V2_newb} U_2 \sim -\frac{|\by|^2}{4}- 2\pi S\nu\,  \by^{T} \mathcal{B}\by + \cdots , \qquad |\by|\to\infty.
\end{gather}
\esub
In Appendix \ref{app:eqnv2}, we state and solve the canonical problem \eqref{eqnmain:V2_new} and obtain the refined behavior 
\bsub\label{eqn:farfieldU2_final}
\begin{equation}
U_2 \sim -\frac{|\by|^2}{4} -2\pi S\nu \by^{T} \mathcal{B} \by + d_{2c} + \mathcal{O}(|\by|^{-2}), \quad \mbox{as} \quad |\by|\to\infty,
\end{equation}
where the constant term is
\begin{equation}
d_{2c} = \frac{\alpha^2 + \beta^2}{4} + 4\pi S\nu\alpha\beta \, \mathcal{B}_{11}.
\end{equation}
\esub
If we incorporate the value of $\mathcal{B}_{11}$ shown in \eqref{eqn:MatrixB_a}, this term reduces to
\begin{equation}\label{eq:d2c_final}
d_{2c} = \frac{\alpha^2 + \beta^2}{4} - \pi S\nu \mbox{Trace}\big(\Q \nabla_{\bx}^2 R(\bxi;\bxi) \big),
\end{equation}
where we have used the identity
\[
  \mbox{Trace}\big(\Q \nabla_{\bx}^2 R(\bxi;\bxi) \big) = \Q_{11} (R_{11} - R_{22}) + 2\Q_{12} R_{12}; \qquad\Q = \begin{bmatrix} \Q_{11} & \phantom{-}\Q_{12}\\ \Q_{12} & -\Q_{11} \end{bmatrix} = -\alpha\beta \begin{bmatrix} \cos2\phi & \phantom{-}\sin2\phi\\ \sin2\phi & -\cos2\phi \end{bmatrix}.
\]

This completes the solution of the inner expansion \eqref{eqn:expInner} the inner problem up to $\bigoh(\eps^2)$. A combination of equations \eqref{eqn:farfieldU0}, \eqref{eqn:farfieldU1} and \eqref{eqn:farfieldU2_final} yields the local behavior 
\begin{align}\label{eqnUlocalFull}
\nonumber u &\sim S\nu \log|\bx-\bxi|  + S \\[4pt]
  &+  \eps^2 \left(S\nu\left[ \frac{(\bx-\bxi)^T \Q(\bx-\bxi) }{|\bx-\bxi|^4} - 2\pi\,\ba \cdot \frac{\M (\bx-\bxi) }{|\bx-\bxi|^2}\right] + d_{2c}  \right) \quad \mbox{as} \quad \bx\to\bxi.
\end{align}
We now return to the outer expansion. 

\underline{Outer region $\bigoh(\eps^1)$}: At this order the problem is given by
\bsub\label{eqn:sing_v1}
\begin{align}
\Delta u_1 &= 0, \qquad \bx \in \Omega\setminus\{\bxi\};\\[5pt]
\nabla u_1 \cdot \hn &= 0, \qquad \bx\in\partial\Omega; \\[5pt]
u_1 &\sim  0, \qquad \bx\to\bxi.
\end{align}
\esub
The unique solution of \eqref{eqn:sing_v1} is $u_1 \equiv 0$.

\underline{Outer region $\bigoh(\eps^2)$}: At this order we have that
\bsub\label{eqn:sing_lv2}
\begin{align}
\label{eqn:sing_v2_a} \Delta u_2 &= 0, \qquad \bx \in \Omega\setminus\{\bxi\};\\[5pt]
\label{eqn:sing_v2_b} \nabla u_2 \cdot \hn &= 0, \qquad \bx\in\partial\Omega;\\[5pt]
\label{eqn:sing_v2_c}  u_2 & \sim S \nu \left[ \frac{(\bx-\bxi)^T \Q(\bx-\bxi)}{|\bx-\bxi|^4} -2\pi \ba \cdot  \frac{\M (\bx-\bxi)}{|\bx-\bxi|^2}\right] + d_{2c} + \cdots, \qquad \bx\to\bxi.
\end{align}
\esub
To express the solution of \eqref{eqn:sing_lv2} in terms of the Green's function, we first notice by direct calculation that
\bsub
\begin{align}
\nabla_{\bxi} \log{|\bx-\bxi|} &= -\frac{\bx-\bxi}{|\bx-\bxi|^2},\\[5pt]
-\frac{1}{2}\mbox{Trace} \big(\Q\nabla_{{\bxi}}^2 \log |\bx-\bxi| \big) &=  \frac{(\bx-\bxi)^T \Q(\bx-\bxi) }{|\bx-\bxi|^4}, \qquad  \Q = \begin{bmatrix} Q_{11} & \phantom{-}Q_{12} \\ Q_{12} & -Q_{11} \end{bmatrix}.
\end{align}
\esub

The solution of \eqref{eqn:sing_lv2} can then be written as
\begin{equation}\label{eqn:s2}
u_2(\bx) =  S\nu \Big [\pi \mbox{Trace} \big( \Q \nabla^2_{{\bxi}} G(\bx;\bxi)\big)  -4\pi^2\ba \cdot \M \nabla_{{\bxi}} G(\bx;\bxi) \Big] + \chi_2,
\end{equation}
where $\chi_2$ is a constant to be determined. In the formulation \eqref{eqn:s2}, the derivatives with respect to the source location $\by=(y_1,y_2)$ are
\[
\nabla_{\bxi} = \begin{bmatrix} \partial_{\xi_1}\\[5pt] \partial_{\xi_2}\end{bmatrix}, \qquad  \nabla^2_{\bxi} = \begin{bmatrix} \partial^2_{\xi_1 \xi_1} & \partial^2_{\xi_1\xi_2}\\[5pt] \partial^2_{\xi_2 \xi_1} &  \partial^2_{\xi_2 \xi_2}  \end{bmatrix}.
\]
This in particular leads to the identity
\begin{equation}\label{eq:TraceForm}
\mbox{Trace} \big( \Q \nabla^2_{\bxi} G(\bx;\by)\big) = Q_{11} (G_{y_1y_1} -G_{y_2y_2}) + 2Q_{12} G_{y_1y_2}.
\end{equation}

To complete the matching, we evaluate \eqref{eqn:s2} as $\bx\to\bxi$. Since $\ba = \nabla_{\bx} R (\bxi;\bxi) = \nabla_{\bxi}R (\bxi;\bxi)$, we calculate
\[
u_2 \sim \chi_2 + S\nu \left[  \frac{(\bx-\bxi)^T \Q (\bx-\bxi)}{|\bx-\bxi|^4} -2\pi \ba \cdot  \frac{\M (\bx-\bxi)}{|\bx-\bxi|^2}  + \pi \mbox{Trace}\big(\Q \nabla_{{\by}}^2 R(\bxi;\bxi) \big) - 4\pi^2 \ba \cdot \M \ba  \right], \qquad \bx\to\bxi.
\]
Matching with \eqref{eqn:sing_v2_c}, we have that
\[
\chi_2 = -S\nu\Big( \pi \mbox{Trace}\big(\Q \nabla_{{\by}}^2 R(\bxi;\bxi) \big) - 4\pi^2 \ba \cdot \M \ba \Big) + d_{2c}.\\
\]
In the case particular to the elliptical trap, we apply \eqref{eqn:M} so that $\M = - \alpha^2 \mathcal{I} + \Q$. In addition, we apply the value of $d_{2c}$ given in \eqref{eq:d2c_final}, together with the symmetry relationship $\nabla_{\bx}^2R(\bxi;\bxi) = \nabla_{\bxi}^2R(\bxi;\bxi)$, and $S\nu = |\Omega|/2\pi$ from \eqref{eqn:leading_order_sol} to further reduce $\chi_2$ to
\begin{align}
\nonumber    \chi_2 &= -S\nu\Big( 2\pi \mbox{Trace}\big(\Q \nabla_{\bxi}^2 R(\bxi;\bxi) \big) + 4\pi^2 \alpha^2 |\ba|^2 - 4\pi^2 \ba \cdot \Q \ba \Big) + \frac{a^2+b^2}{8}\\
\label{eqn:chi2_final}   &= -|\Omega|\Big(  \mbox{Trace}\big(\Q \nabla_{\bxi}^2 R(\bxi;\bxi) \big) + 2\pi \alpha^2 |\ba|^2 - 2\pi \ba \cdot \Q \ba \Big) + \frac{a^2+b^2}{8}
\end{align}

In summary, the solution of \eqref{eqn:MFPTIntro} admits the expansion $u = \frac{1}{D}[u_0 + \eps^2 u_2+\cdots]$ where the terms are
\bsub
\begin{align}
u_0(\bx;\bxi) &= -|\Omega| \Big[ G(\bx;\bxi) - R(\bxi;\bxi) \Big] + \frac{|\Omega|}{2\pi \nu},\\[4pt]
u_2(\bx;\bxi) &=  |\Omega|\Big [ \frac12 \mbox{Trace} \big( \Q \nabla^2_{\bxi} G(\bx;\bxi)\big) -2\pi\ba \cdot \M \nabla_{{\bxi}} G(\bx;\bxi) \Big] + \chi_2.
\end{align}
\esub
where $\chi_2$ is the constant given in \eqref{eqn:chi2_final}.
\subsection{Calculation of the global MFPT}
In this section we calculate the GMFPT defined as 
\[
\tau = \frac{1}{|\Omega\setminus\Omega_{\eps}|}\int_{\Omega\setminus\Omega_{\eps}}u \, d\bx =  \frac{1}{D|\Omega\setminus\Omega_{\eps}|}\int_{\Omega\setminus\Omega_{\eps}} u_0 + \eps^2 u_2 \, d\bx.
\]
The challenge as before is to determine the correct expansion accurate to $\bigoh(\eps^2)$ by accounting for contributions from the inner expansion. We first decompose the region of integration  $\Omega\setminus \Omega_{\eps} := (\Omega \setminus B_{\delta})\cup(B_{\delta}\setminus \Omega_{\eps})$ for the disk $B_{\delta} = \{ \bx\in\mathbb{R}^2 \ | \ |\bx-\bxi|<\delta\}$ and then apply the limit $\delta\to0$. The integral becomes
\[
\int_{\Omega\setminus \Omega_{\eps}} u \, d\bx = \underbrace{\int_{B_{\delta}\setminus \Omega_\eps} u \, d\bx}_{I_1} + \underbrace{\int_{\Omega \setminus B_{\delta}} u \, d\bx}_{I_2}.
\]
The integral $I_1$ is evaluated by first transforming to the coordinate $\by = e^{-i\phi}(\bx-\bxi)/\eps$ followed by $\by = \alpha \bz + \beta/\bz$ to find that 
\begin{equation}
I_1 = \eps^2 \int_{\substack{\by\in\mathbb{R}^2\setminus\A\\ |\by|<\delta/\eps}} U\, d\by = \frac{\eps^2}{D} \int_{|\bz|=1}^{|\bz| = \delta/\eps} (\underbrace{U_0}_{I_{10}} + \underbrace{ \eps U_1}_{{I_{11}}} + \underbrace{ \eps^2 U_2}_{I_{12}} + \cdots )|J|d\bz.
\end{equation}
From equation \eqref{eqn:outerS}, we have that $U_0(\bz)= S\nu \log |\bz|$ so that for $\bz = r e^{i \theta}$,
\begin{align}
\nonumber I_{10} &= S\nu\eps^2\int_{\theta = 0}^{2\pi}\int_{r=1}^{\frac{\delta}{\eps\alpha} } \log r \left[ \alpha^2- \frac{2\alpha\beta}{r^2}\cos 2\theta + \frac{\beta^2}{r^4} \right] r dr d\theta = 2\pi S\nu \eps^2 \int_{r=1}^{r=\frac{\delta}{\alpha\eps}}\log r \left[ \alpha^2 r + \frac{\beta^2}{r^3} \right] dr \\[4pt]
\label{eq:Int10} & = 2\pi S\nu \eps^2 \Big[\alpha^2 \Big(\frac{r^2}{2}\log r -\frac{r^2}{4}\Big) - \beta^2 \Big( \frac{1+ 2\log r}{4 r^2} \Big)  \Big]_{r=1}^{r=\frac{\delta}{\eps\alpha}} \approx \frac{S\nu \pi}{2} \left[  2\delta^2\log \frac{\delta}{\eps\alpha}-\delta^2 + \eps^2 (\alpha^2 + \beta^2)\right] .
\end{align}
Following on to the term $I_{11}$, we apply from \eqref{eqn:U1} that $U_1 = -2\pi S \nu \ba \cdot e^{-i\phi} v_{1c} $ where $v_{1c}(\bz) = (\bz - \bz/|\bz|^2)$. This leads to $I_{11}  = \eps^3 \int_{|\bz|=1}^{|\bz|= \delta/\eps} U_1 |J| d\bz = 0$. The contribution from $I_{22} = \bigoh(\eps^4)$.

The contribution from the outer region is now
\begin{equation}\label{eqn:I2}
I_2 = \int_{\Omega \setminus B_{\delta}} u \, d\bx = \frac{1}{D} \underbrace{\int_{\Omega \setminus B_{\delta}} u_0 \, d\bx}_{I_{20}} + \frac{\eps^2}{D}  \underbrace{\int_{\Omega \setminus B_{\delta}} u_2 \, d\bx}_{I_{22}}.
\end{equation}
The first term is calculated with $u_0 = -2\pi S\nu G(\bx;\bxi) + \tau_0$ with $\tau_0 = S(1+ 2\pi \nu R(\bxi;\bxi))$ to determine
\begin{align}
\nonumber I_{20} &= \int_{\Omega \setminus B_{\delta} } u_0\, d \bx = \int_{\Omega} u_0\,d\bx -  \int_{B_{\delta}} u_0\, d\bx \\[4pt]
\nonumber {} &=  |\Omega| \tau_0 - \int_{B_{\delta} } [S\nu  \log |\bx-\bxi| +S]\, d\bx -2\pi S\nu  \int_{B_{\delta} } (R(\bx;\bxi)- R(\bxi;\bxi) )\, d\bx\\[4pt]
\label{eqn:Int20} {} &= |\Omega| \tau_0 - S\nu \pi \delta^2\Big(\log \delta - \frac12\Big) - S\pi \delta^2 + \bigoh(\delta^4).
\end{align}
Following on to the second term in \eqref{eqn:I2}, we calculate that
\begin{align}
\nonumber I_{22} = \int_{\Omega \setminus B_{\delta} } u_2 \, d\bx = \int_{\Omega} u_2\,d\bx -  \int_{B_{\delta}} u_2\, d\bx.
\end{align}
From \eqref{eqn:sing_v2_b}, the contributions of the second term have vanishing average, hence
\[
\int_{B_{\delta}} u_2\, d\bx = \pi \delta^2 \, \frac{(a^2+b^2)}{8}.
\]
Combining \eqref{eq:Int10} and \eqref{eqn:Int20}, we have that
\begin{align}
\nonumber \int_{\Omega\setminus\Omega_{\eps}} u_0 + \eps^2 u_2 \, d\bx &=  |\Omega|\tau_0 + \eps^2 \Big( \frac{S\nu\pi}{2} (\alpha^2 + \beta^2) + \int_{\Omega} u_2 \, d\bx \Big) + \bigoh(\eps^3),\\
&= |\Omega| \tau_0 + \eps^2 |\Omega| \Big( \frac{a^2+b^2}{8} + \chi_2 +  \int_{\Omega} \frac{1}{2} \mbox{Trace} \big( \Q \nabla^2_{{\bxi}} G(\bx;\bxi)\big)  -2\pi\ba \cdot \M \nabla_{{\bxi}} G(\bx;\bxi) d\bx \Big)+ \bigoh(\eps^3). \label{eqn:int2}
\end{align}
Applying the identity \eqref{eq:TraceForm}, we note that the integral terms in \eqref{eqn:int2} vanish since
\begin{equation}
\int_{\Omega} (G_{\xi_1\xi_1} - G_{\xi_2\xi_2}) d \bx = 0, \quad \int_{\Omega} G_{\xi_1\xi_2}  d \bx = 0, \quad \int_{\Omega} G_{\xi_1}  d \bx = 0, \quad \int_{\Omega} G_{\xi_2}  d \bx =0.
\end{equation}
Hence we obtain the final expression which incorporates the local and global contributions to the GMFPT,
\begin{equation}
\int_{\Omega\setminus\Omega_{\eps}} u_0 + \eps^2 u_2 \, d\bx = |\Omega| \tau_0 + \eps^2 |\Omega| \Big( \frac{a^2+b^2}{8} + \chi_2 \Big).
\end{equation}

Recalling the value of $\chi_2$ given in \eqref{eqn:chi2_final}, equation \eqref{eqn:int2} then reduces to 
\begin{equation}\label{eqn:deftau}
\int_{\Omega\setminus\Omega_{\eps}} u_0 + \eps^2 u_2 \, d\bx = |\Omega| \tau_0 + \eps^2 |\Omega| \left( \frac{a^2+b^2}{4} - |\Omega|\Big( \mbox{Trace}\big(\Q \nabla_{\bxi}^2 R(\bxi;\bxi) \big) + \pi \frac{(a+b)^2}{2} |\ba|^2  - 2\pi \ba \cdot \Q \ba \Big) \right).
\end{equation}
Finally, using $|\Omega\setminus\Omega_{\eps}| = |\Omega| -\pi \eps^2 ab$, we have that
\begin{align}
\nonumber  \tau &=  \frac{1}{D|\Omega\setminus\Omega_{\eps}|}\int_{\Omega\setminus\Omega_{\eps}} u_0 + \eps^2 u_2 d\bx = \frac{1}{D|\Omega|}\left[ 1 + \eps^2 \frac{\pi a b}{|\Omega|} \right]\int_{\Omega\setminus\Omega_{\eps}} u_0 + \eps^2 u_2 \, d\bx\\
    &= \frac{\tau_0}{D} + \frac{\eps^2}{D} \left( \frac{\pi ab}{|\Omega|}\tau_0 +  \frac{a^2+b^2}{4} - |\Omega|\Big( \mbox{Trace}\big(\Q \nabla_{\bxi}^2 R(\bxi;\bxi) \big) + \pi \frac{(a+b)^2}{2} |\ba|^2  - 2\pi \ba \cdot \Q \ba \Big) \right). \label{eqn:tau_final}
\end{align}

\begin{figure}[htbp]
\centering
\subfigure[Variation in $u_2(\bx)=\eps^{-2}(u(\bx)-u_0(\bx))$ with $\phi$.]{\includegraphics[width = 0.415\textwidth]{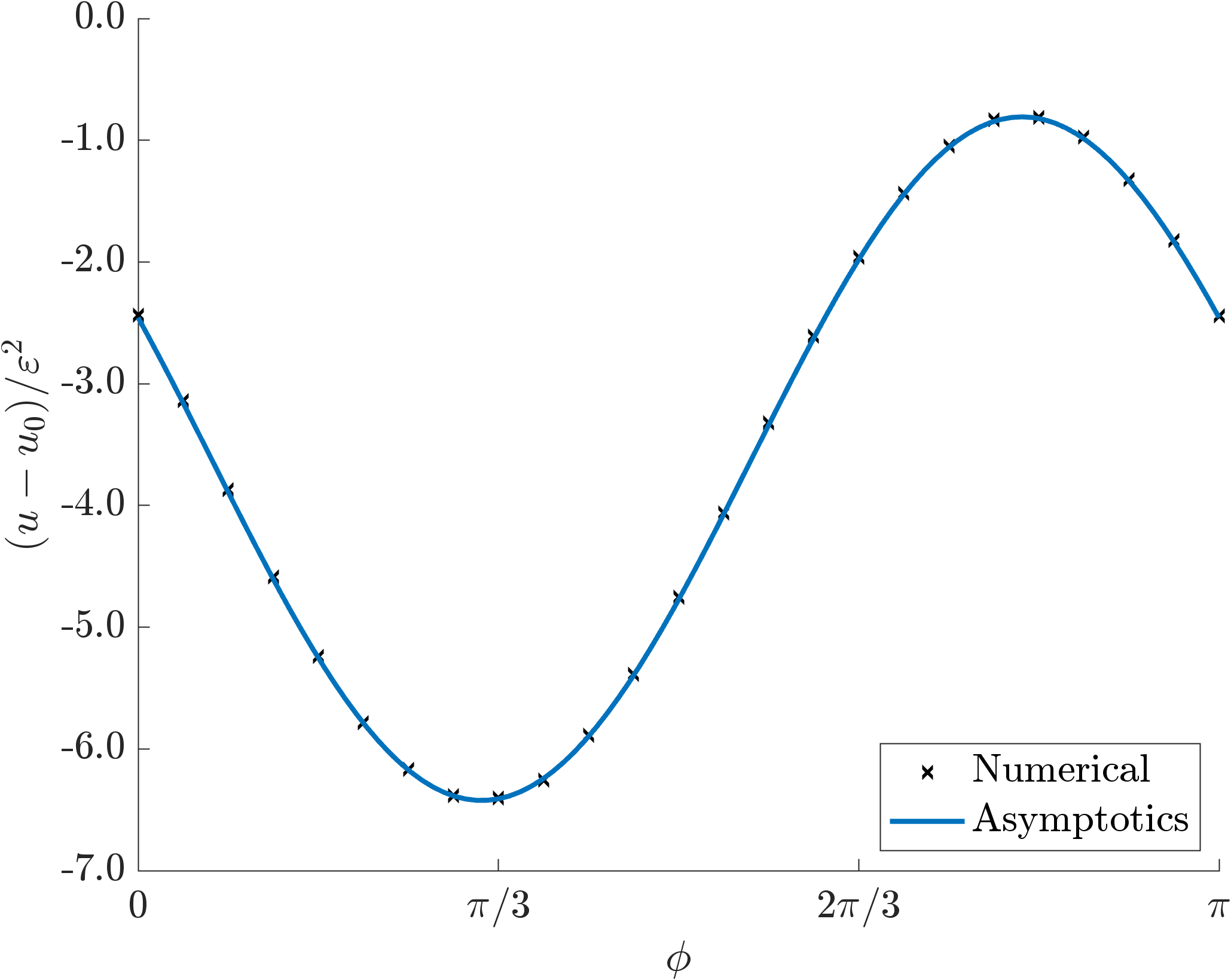}\label{fig_convergence_a}}\hspace{1cm}
\subfigure[Variation in $\tau_2 = \eps^{-2}(\tau - \tau_0)$ with $\phi$.]{\includegraphics[width = 0.415\textwidth]{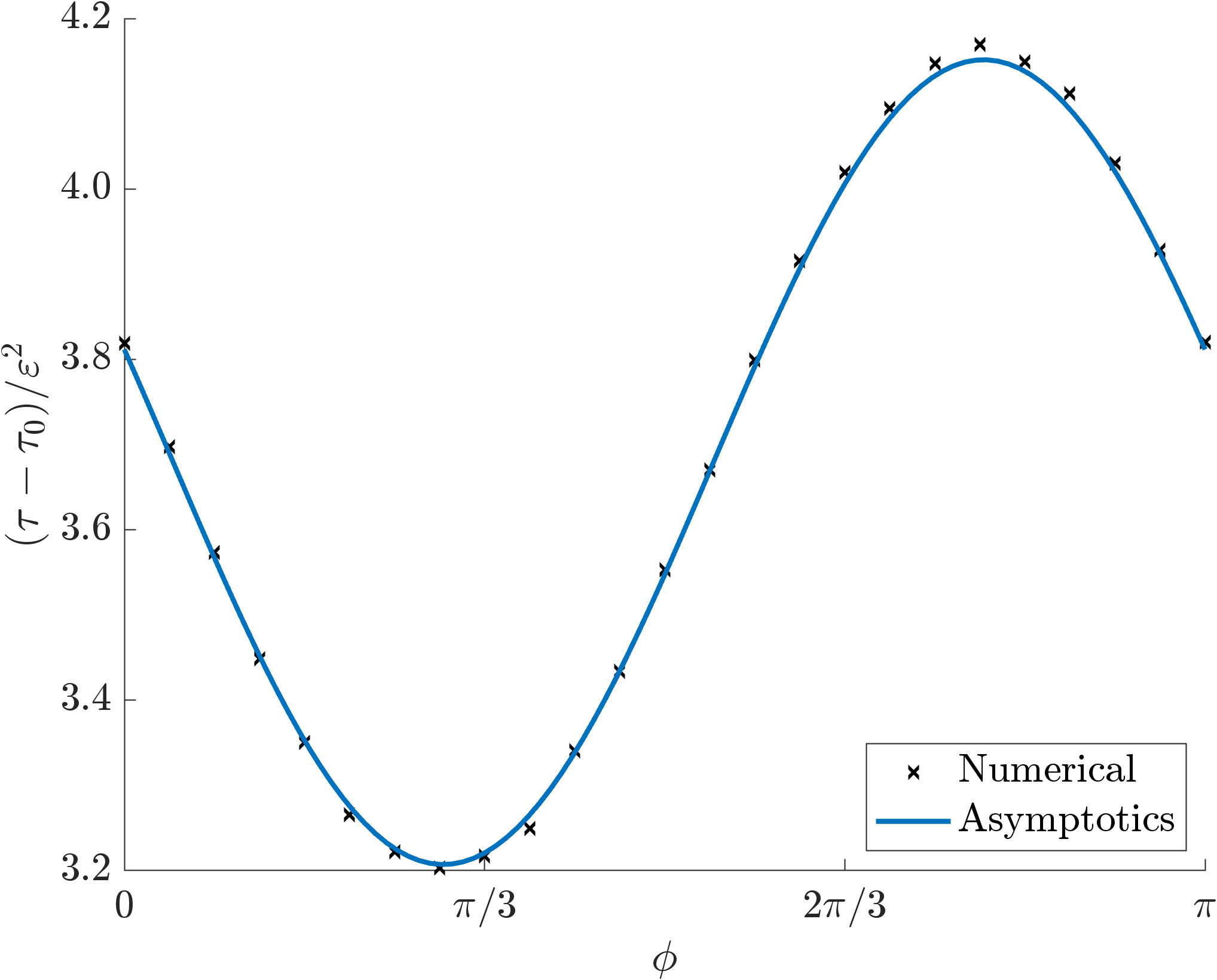}\label{fig_convergence_b}}\\
\subfigure[Convergence of $u(\bx)$.]{\includegraphics[width = 0.415\textwidth]{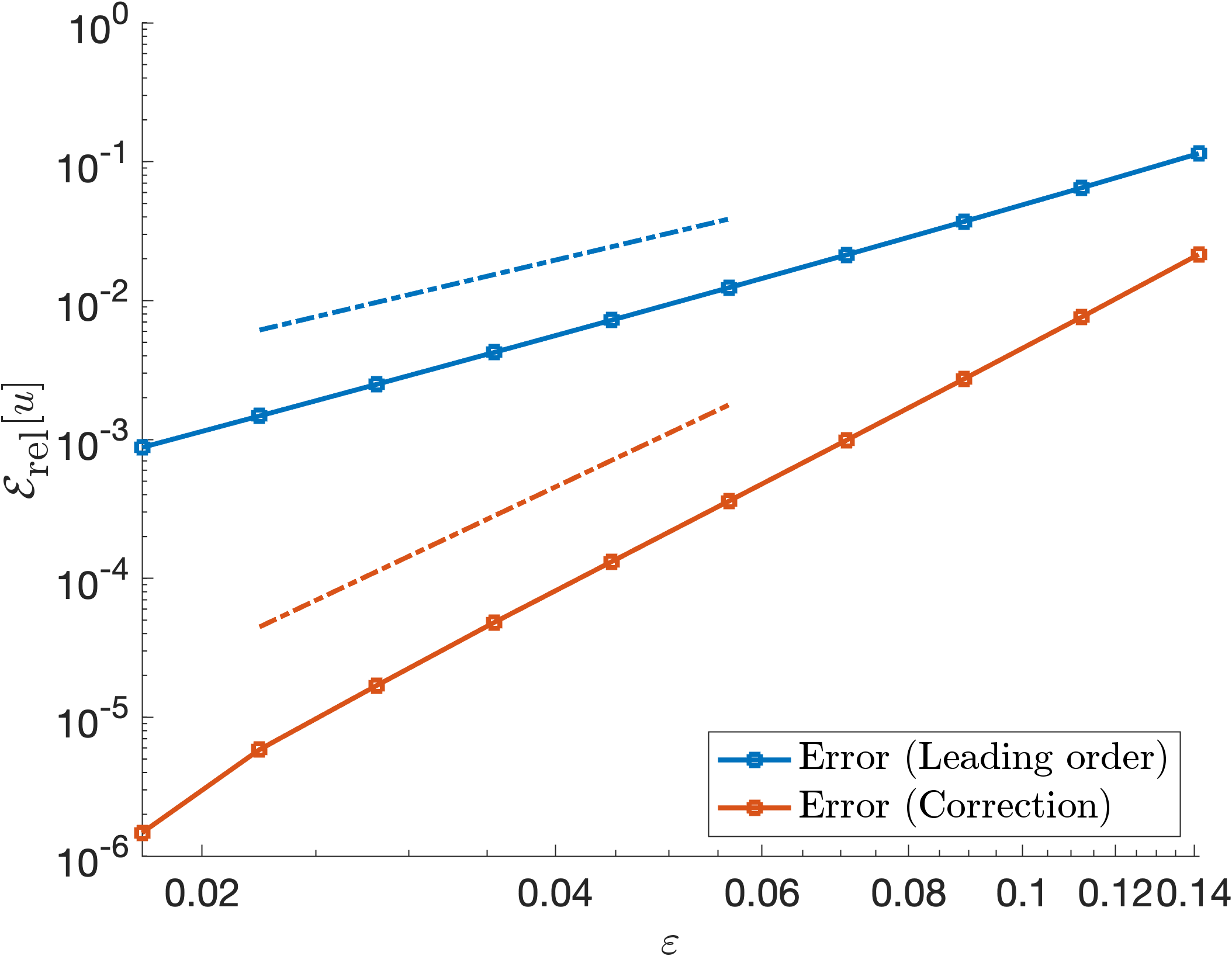}\label{fig_convergence_c}}\hspace{1cm}
\subfigure[Convergence of $\tau$.]{\includegraphics[width = 0.415\textwidth]{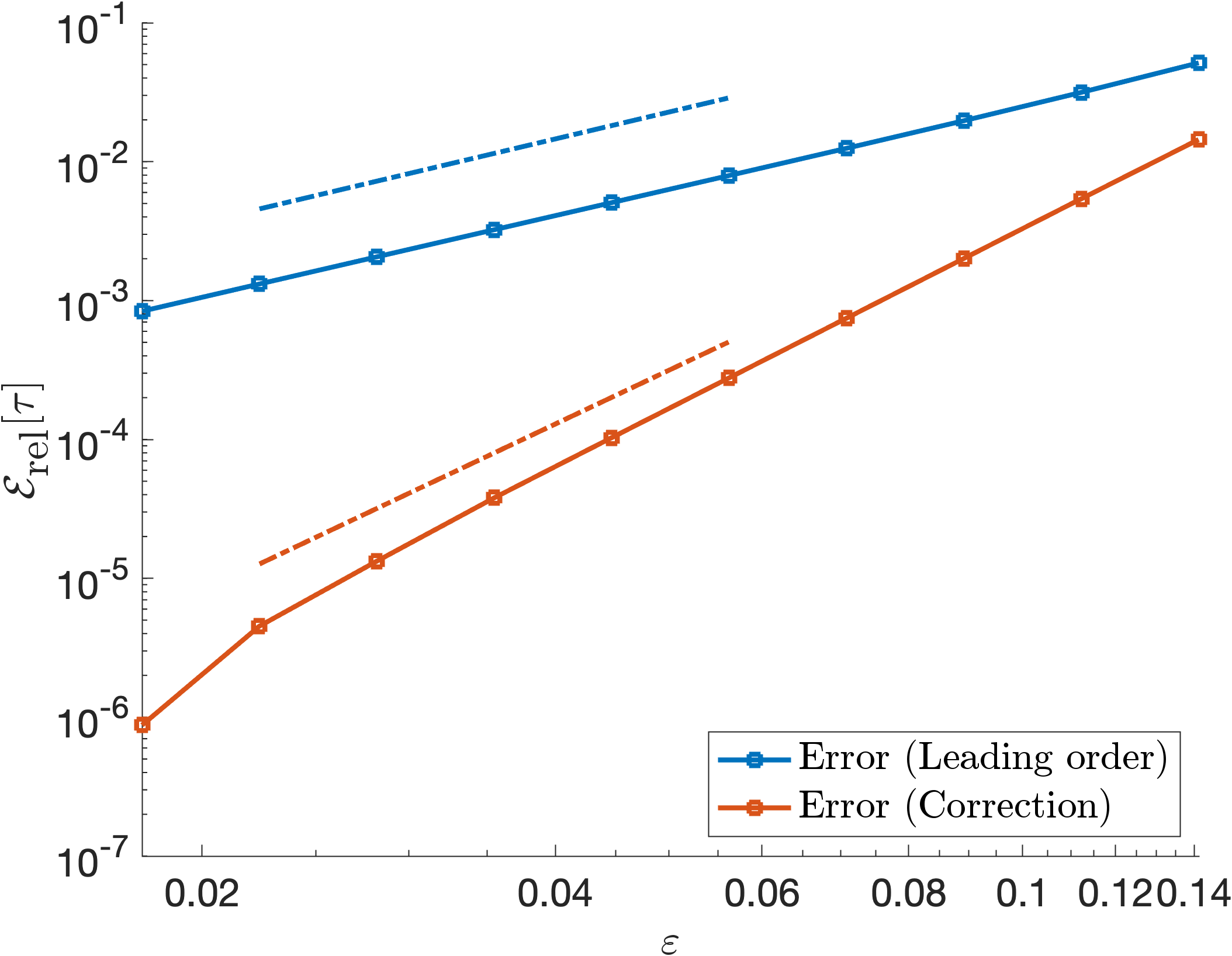}\label{fig_convergence_d}}
\caption{Convergence of the asymptotic approximation \eqref{eq:introExpansion} in the disk case with a trap centered at $\bxi = (0.3,0.4)$. Panel (a): Agreement between the solution correction $u_2(\bx) = \eps^{-2}(u(\bx)-u_0(\bx))$ for $\bx = (-0.2,-0.4)$. Panel (b): The GMFPT correction $\tau_2= \eps^{-2}(\tau-\tau_0)$ from numerical and asymptotic approximations for $\eps = 0.03$ as orientation $\phi$ varies. Panels (c-d): Convergence as $\eps \to 0$ of the relative errors between numerical and asymptotic approximations (leading and correction) of $u(\bx)$ for $\bx = (-0.2,-0.4)$ (c), and $\tau$ with fixed $\phi = \pi/6$ (d). Straight lines are of slope $2$ (blue) and $4$ (red) indicating convergence rates. Domain schematic shown in Fig.~\ref{fig:optim_g_a}. \label{fig_convergence}}
\end{figure}

\section{Results}\label{sec:results}

In this section, we demonstrate the validity of the expansion and investigate how trap orientation modulates the GMFPT. In terms of minimizing the GMFPT, we remark that trap location is felt in the leading order through the term $R(\bxi;\bxi)$. Hence, when considering the role of orientation, we primarily focus on the correction terms $u_2(\bx)$ and $\tau_2$. A globally optimizing configuration of the GMFPT would first locate the trap at the critical points of $\tau_0$ followed by choosing the orientation that optimizes $\tau_2$.

\begin{exmp}
    Convergence in the disk geometry with a single elliptical trap.
\end{exmp}

\begin{figure}[htbp]
    \centering
    \subfigure[MFPT $u(\bx)$ with elliptical trap at $\bxi = (0.3,0.4)$.]{\includegraphics[width = 0.415\textwidth]{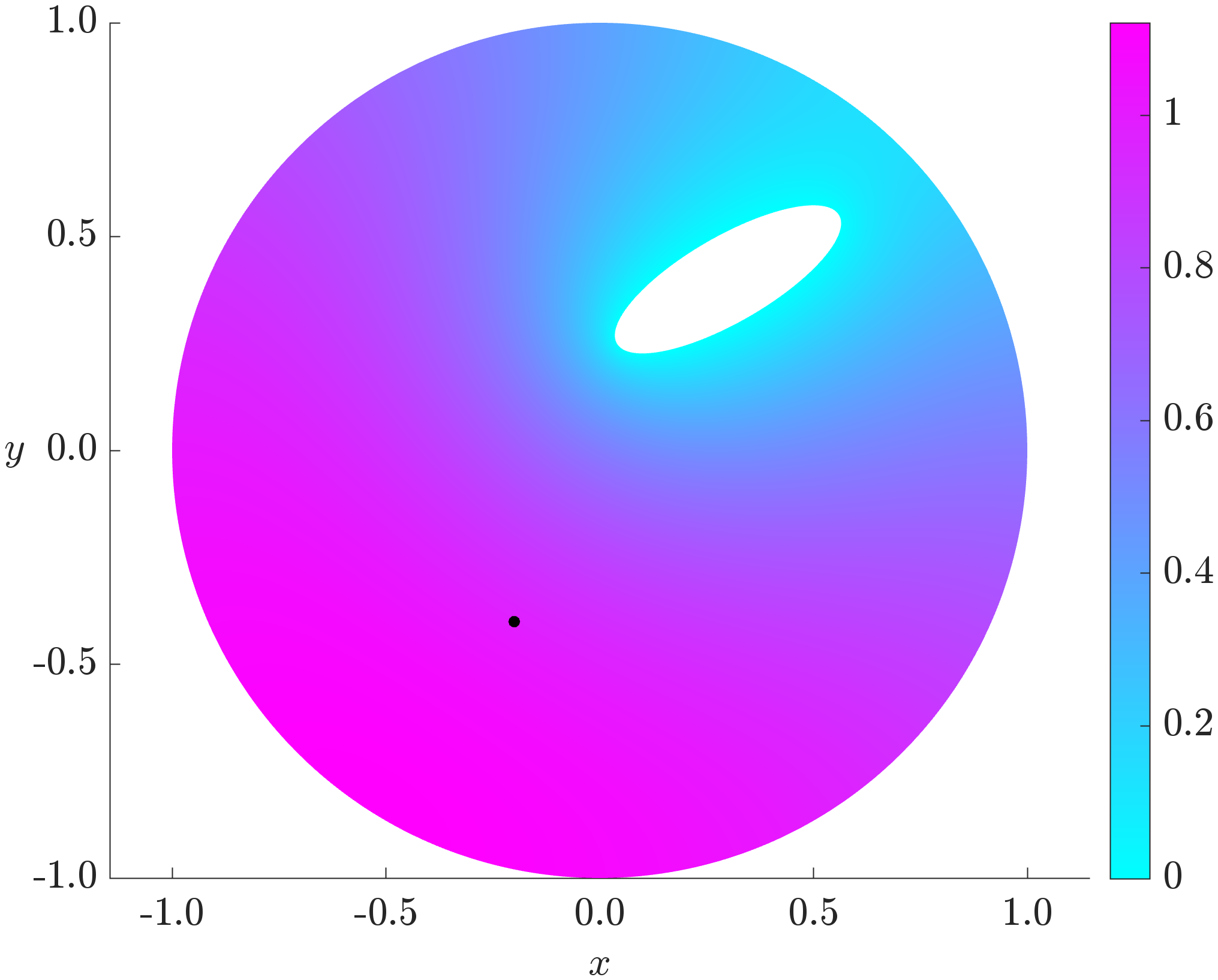}\label{fig:optim_g_a}}\hspace{1cm}
    \subfigure[Orientations that minimize the GMFPT correction $\tau_2$.]{\includegraphics[width = 0.5\textwidth]{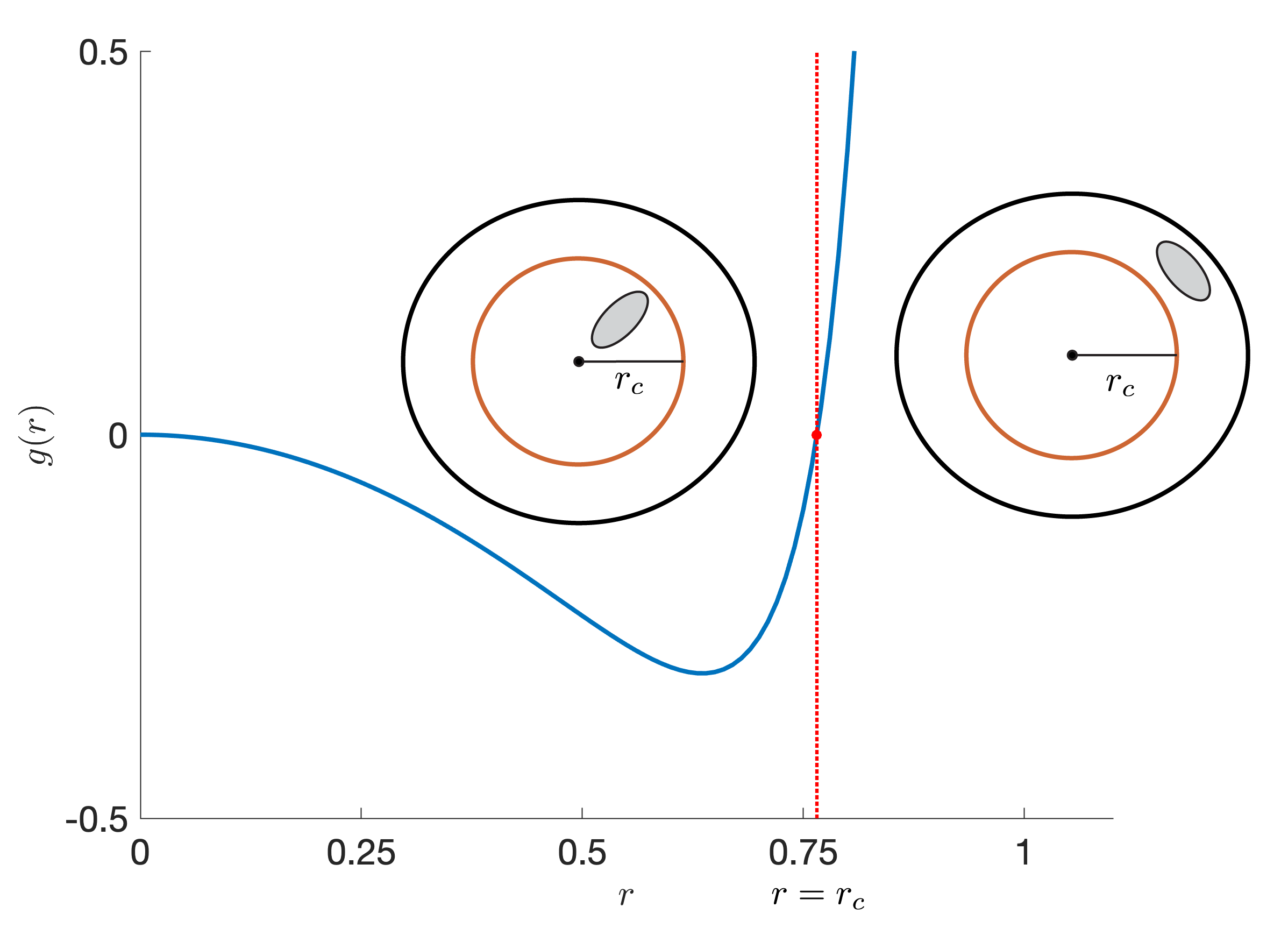}\label{fig:optim_g_b}}
    \caption{Minimization of the GMFPT correction in the disk with a single elliptical trap. Panel (a): Domain with a single elliptical trap at $\bxi = (0.3,0.4)$, axes $\eps(a,b) = \eps(3,1)$ and orientation $\phi = \pi/6$. The highlighted point (black dot) is $\bx=(-0.2,-0.4)$. Panel (b): The function $g(r)$ and the critical radius $r=r_c$. For $r<r_c$, the GMFPT correction is minimized when the ellipse has major axis pointed towards the center of the disk.\label{fig:optim_g}}
\end{figure}

In this case where the enclosing geometry $\Omega$ is the unit disk, we have exact formula for the Green's function and we can obtain the necessary derivatives. From the equations \eqref{eq:AppGreensDisk}, we calculate that
\begin{gather}\label{sol_disk}
    \nonumber  \ba = \nabla_{\bxi}R(\bxi;\bxi) = \frac{1}{2\pi}\frac{ 2 - |\bxi|^2}{1-|\bxi|^2} \bxi; \quad 
     \mbox{Trace}\big(\Q \nabla^2_{\bxi} R(\bxi;\bxi) \big)  = \frac{1}{\pi} \frac{\bxi \cdot \Q \bxi}{(1-|\bxi|^2)^2}\, ; \quad 
     \ba \cdot \Q \ba  =\frac{1}{4\pi^2} \frac{(2-|\bxi|^2)^2}{(1-|\bxi|^2)^2}\, \bxi \cdot \Q \bxi;\\[5pt]
     S\nu = \frac{|\Omega| }{2\pi } = \frac{1 }{2}; \qquad \tau_0 = \frac{|\Omega|}{2\pi \nu} \Big[1 + 2\pi\nu R(\bxi;\bxi) \Big] = \frac{1}{2} \Big[ \frac{1}{\nu} - \log(1-|\bxi|^2) + |\bxi|^2 -\frac34 \Big].
\end{gather}
Applying \eqref{sol_disk} to \eqref{eqn:tau_final} together with $|\Omega| = \pi$ and $D=1$, yields the GMFPT
\begin{align}
\nonumber \tau & = \tau_0 + \eps^2 \left[ ab\tau_0 + \frac{a^2+b^2}{4} - \frac{(a+b)^2}{8} \left(\frac{2-|\bxi|^2}{1-|\bxi|^2}\right)^2 |\bxi|^2 - \frac{2- (2- |\bxi|^2)^2}{2(1-|\bxi|^2)^2} \bxi \cdot \Q \bxi \right]\\ 
& = \tau_0 + \eps^2  \left[ab\tau_0 + \frac{a^2+b^2}{4} - \frac{(a+b)^2}{8} \left(\frac{2-|\bxi|^2}{1-|\bxi|^2}\right)^2 |\bxi|^2 + \frac{a^2-b^2}{4} \frac{2 - (2- |\bxi|^2)^2}{2(1-|\bxi|^2)^2} \begin{bmatrix} \xi_1^2 - \xi_2^2\\ 2\xi_1 \xi_2 \end{bmatrix} \cdot \begin{bmatrix} \cos 2\phi \\ \sin 2\phi \end{bmatrix}\right]. \label{eqn:tauDisk} 
\end{align}
In Fig.~\ref{fig_convergence}, we show agreement between finite element solutions of \eqref{eqn:MFPTIntro} and the asymptotic result \eqref{eqn:tauDisk}. The common variables in this validation are the disk geometry, the trap center $\bxi = (0.3,0.4)$ and the semi-axes dimensions $(a,b)=(3,1)$. In Figs.~\ref{fig_convergence_a}-\ref{fig_convergence_b} we show agreement of the correction terms $u_2 = (u(\bx)- u_0(\bx))/\eps^2$  at the point $\bx = (-0.2,-0.4)$ and $\tau_2 = (\tau - \tau_0)/\eps^2$ as the trap orientation $\phi$ varies with $ \eps = 0.03$ fixed. A schematic of the domain is shown in Fig.~\ref{fig:optim_g_a}.

To directly confirm the convergence of the asymptotic expansion, we calculate a sequence of relative errors
\begin{equation}\label{eqn:rel_error}
\mathcal{E}_{\textrm{rel}}[z] = \left|\frac{z_{\text{approx}}(\eps) - z_{\text{true}}}{z_{\text{true}}}\right|.
\end{equation}
at a range of $\eps$ values. In \eqref{eqn:rel_error} the true value is calculated from finite element simulations of \eqref{eqn:MFPTIntro} and the approximate values come from taking one or two terms in the expansion \eqref{eq:introExpansion_a}. In Figs.~\ref{fig_convergence_c}-\ref{fig_convergence_d}, we show that the rate of convergence in $\mathcal{E}_{\textrm{rel}}$ as $\eps\to0^{+}$ is in agreement with the principal result \eqref{eq:introExpansion}, and in particular, the two term approximation \eqref{eqn:tauDisk} of the GMFPT has error $\bigoh(\eps^4)$.

\begin{exmp}
    Optimization of trap orientation for the GMFPT in the unit disk domain.
\end{exmp}
To study the optimizing trap orientations, we set $\bxi = r e^{i\theta}$ in \eqref{eqn:tauDisk} which reduces the GMFPT to
\begin{equation}\label{eqn:tauDisk_r}
\tau \sim \tau_0 + \eps^2 \left[ ab\tau_0 +\frac{a^2+b^2}{4} - \frac{(a+b)^2}{8} r^2\left(\frac{2-r^2}{1-r^2}\right)^2 + \frac{a^2-b^2}{4} g(r) \cos2(\theta-\phi) \right], \quad g(r)= \frac{2 - (2- r^2)^2}{2(1-r^2)^2}r^2.
\end{equation}

The function $g(r)$ has the following (see Fig.~\ref{fig:optim_g_b}) simple properties, 
\[
g(r) \to 0^{-},\quad \text{as}\quad  r \to 0^{+}; \qquad g(r)\to +\infty, \quad \text{as}\quad  r \to 1^{-}; \qquad g(r_c) = 0, \quad r_c = \sqrt{2-\sqrt{2}},
\]
and hence we can conclude that the correction $\tau_2$ of the GMFPT is minimized when the semi-major axis of the trap is aligned in the radial direction, i.e.~$\phi = \theta$, provided $g(r)<0$ or $(2-r^2)^2>2$. The GMFPT minimizing configuration is therefore
\begin{equation}
    \phi = \left\{
    \begin{array}{rl}
        \theta, & r \in (0,r_c),\\[5pt]
        \theta+\frac{\pi}{2}, & r \in (r_c,1);
    \end{array} \right. \qquad r_c = \sqrt{2-\sqrt{2}} \approx 0.7654.
\end{equation}
In Fig.~\ref{fig:optim_g_b}, we plot the function $g(r)$ together with corresponding trap orientations that minimize the GMFPT. We also remark that the trap center which minimizes the leading order term $\tau_0$ in GMFPT is $\bxi = (0,0)$. For a trap centered at this location, there is no contribution ($g(0)=0$) to the GMFPT from the trap orientation as expected by symmetry considerations.

\begin{figure}[htbp]
  \centering
\includegraphics[width=0.95\textwidth]{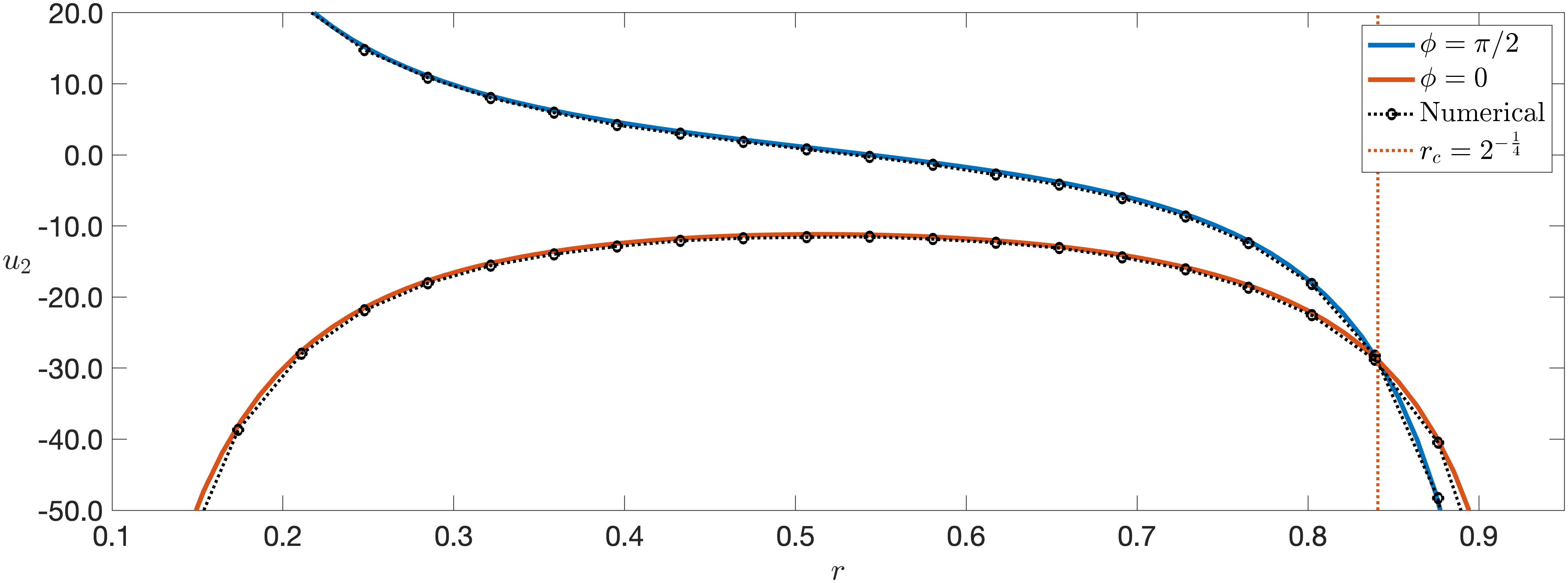} 
\caption{The effects of trap orientation on the MFPT staring at the center of a disk. The correction $u_2 = \eps^{-2}(u-u_0)$ from \eqref{eqn:MFPT_u2} with a single trap of extent $\eps = 0.01$, semi-major axes $(a,b)=(3,1)$ centered at $\bxi=(r,0)$. Curves shown for orientations $\phi =\pi/2$ and $\phi = 0$. \label{fig:CenterEscape}}
\end{figure}

\begin{exmp}
The MFPT from the origin to an elliptical trap in a disk.
\end{exmp}
In this example, we calculate the MFPT for a particle in the disk domain initially at $\bx=0$ to arrive at an elliptical trap centered at $\bxi = re^{i\theta}$. From the main result \eqref{eqn:intro_u2}, we have the correction term
\[
u_2(\bx;\bxi) =  |\Omega|\Big [ \frac12 \mbox{Trace} \big( \Q \nabla^2_{\bxi} G(\bx;\bxi)\big) -2\pi\nabla_{\bxi}R(\bxi;\bxi) \cdot \M \nabla_{{\bxi}} G(\bx;\bxi) \Big] + \chi_2,
\]
which describes the role that trap orientation plays in the capture rate. In Appendix \ref{sec:greens_appendix}, we calculate the relevant terms to be
\bsub
\begin{align}
   \text{Trace}(\mathcal{Q} \nabla_{\bxi}^2G(\textbf{0};\bxi)) &= \frac{1}{\pi}\frac{\bxi \cdot \mathcal{Q} \bxi}{|\bxi|^4};\\
     \nabla_{\bxi}R(\bxi;\bxi)\cdot \mathcal{M}\nabla_{\bxi}G(\textbf{0};\bxi) &= -\frac{1}{4\pi^2}\frac{2-|\bxi|^2}{|\bxi|^2} \left[ \bxi\cdot \mathcal{Q} \bxi - \frac{(a+b)^2}{4} |\bxi|^2 \right].
\end{align}
We then have that 
\begin{equation}\label{eqn:MFPT_u2}
    u_2(\textbf{0}) = \frac{a^2+b^2}{8} - \frac{(a+b)^2}{8} \frac{(2-r^2)}{(1-r^2)^2} +\frac{a^2-b^2}{8} \frac{2r^4-1}{r^2(1-r^2)^2} \cos\big( 2(\phi-\theta)\big).
\end{equation}
\esub
Favorable agreement between this result and numerical simulations is shown in Fig.~\ref{fig:CenterEscape}. We conclude that $u_2(\textbf{0})$ is minimized by orientating the trap in the radial direction ($\theta = \phi$) for $2r^4-1<0$ or $0<r< 2^{-\frac14}\approx 0.8409$. For $2^{-\frac14}<r<1$, $u_2(\textbf{0})$ is minimized by orientating the semi-major axis of the trap parallel to the boundary.

\begin{exmp}
    Here we consider a general domain $\Omega$ with a single elliptical trap centered at $\bx=\bxi$ of semi-major and semi-minor axes $\eps a$ and $\eps b$ respectively and orientation $\phi$ with respect to the horizontal.
\end{exmp}
The GMFPT \eqref{eqn:deftau} reduces to the form
\begin{equation}\label{eqn:GlobalTauSimp}
\tau = \tau_0 + \eps^2 \left[\frac{\pi ab}{|\Omega|}\tau_0 + \frac{a^2+b^2}{4} - \pi |\Omega| \frac{(a+b)^2}{2} (R_{\xi_1}^2 + R_{\xi_2}^2) + |\Omega| \frac{a^2-b^2}{4} \bp \cdot \begin{bmatrix}\cos 2\phi \\ \sin 2\phi \end{bmatrix} \right].
\end{equation}
The vector $\bp$ is given as
\begin{equation}\label{eqn:example_p}
\bp = \begin{bmatrix} 
R_{\xi_1\xi_1}-R_{\xi_2\xi_2} - 2\pi (R_{\xi_1}^2 -R_{\xi_2}^2) \\[4pt]
2R_{\xi_1\xi_2} - 4\pi R_{\xi_1} R_{\xi_2}
\end{bmatrix} = \begin{bmatrix} 
R_{\xi_1\xi_1}-R_{\xi_2\xi_2} \\[4pt]
2R_{\xi_1\xi_2} 
\end{bmatrix} - 2\pi\begin{bmatrix} 
   R_{\xi_1}^2 -R_{\xi_2}^2 \\[4pt]
  2 R_{\xi_1} R_{\xi_2}
\end{bmatrix}.
\end{equation}
For some regular domains such as ellipses and rectangles, the function $R(\bx;\bxi)$ is available in the form of rapidly convergent series (see Appendices  \ref{sec:GreensRect} and \ref{sec:GreensEllipse}). The necessary derivatives $R_{\xi_1\xi_1},R_{\xi_2\xi_2},R_{\xi_1\xi_2},R_{\xi_1}$ and $R_{\xi_2}$ in \eqref{eqn:example_p} can then be calculated by finite differences. We remark that if the trap center is chosen to optimize $\tau_0$, then we have $R_{\xi_1} = R_{\xi_2} = 0$ and so 
\begin{equation}\label{eqn:GlobalTauSimpReduced}
\tau = \tau_0 + \eps^2 \left[\frac{\pi ab}{|\Omega|}\tau_0 + \frac{a^2+b^2}{4} + |\Omega| \frac{a^2-b^2}{4} \begin{bmatrix} 
R_{\xi_1\xi_1}-R_{\xi_2\xi_2} \\[4pt]
2R_{\xi_1\xi_2}.
\end{bmatrix} \cdot \begin{bmatrix}\cos 2\phi \\ \sin 2\phi \end{bmatrix} \right].
\end{equation}

\begin{figure}[htbp]
\centering
\subfigure[$L = 1$, $d = 1$.]{\includegraphics[height = 1.8in]{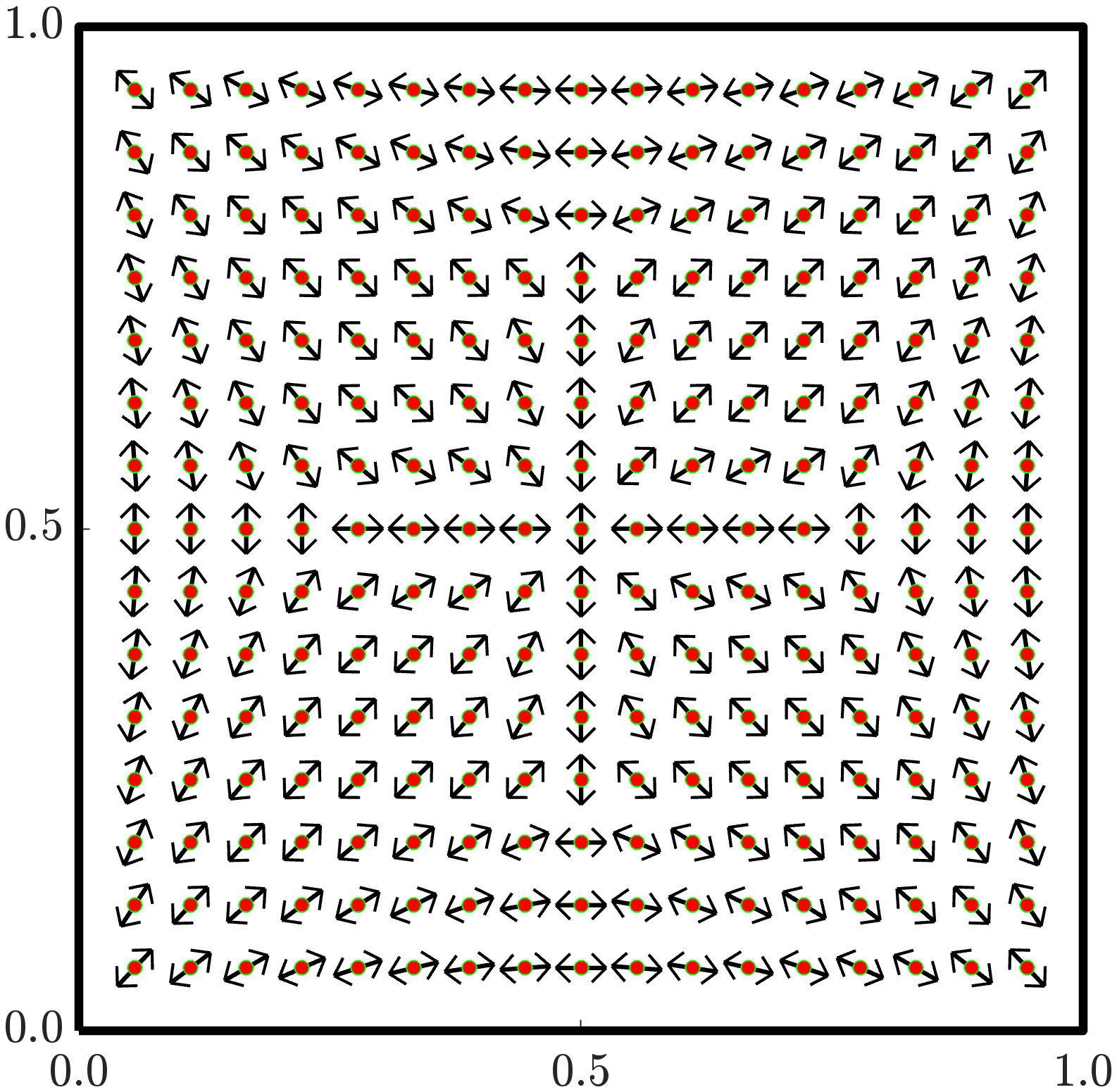}}\hspace{0.5cm}
\subfigure[$L = 1.01$, $d = 1$.]{\includegraphics[height = 1.8in]{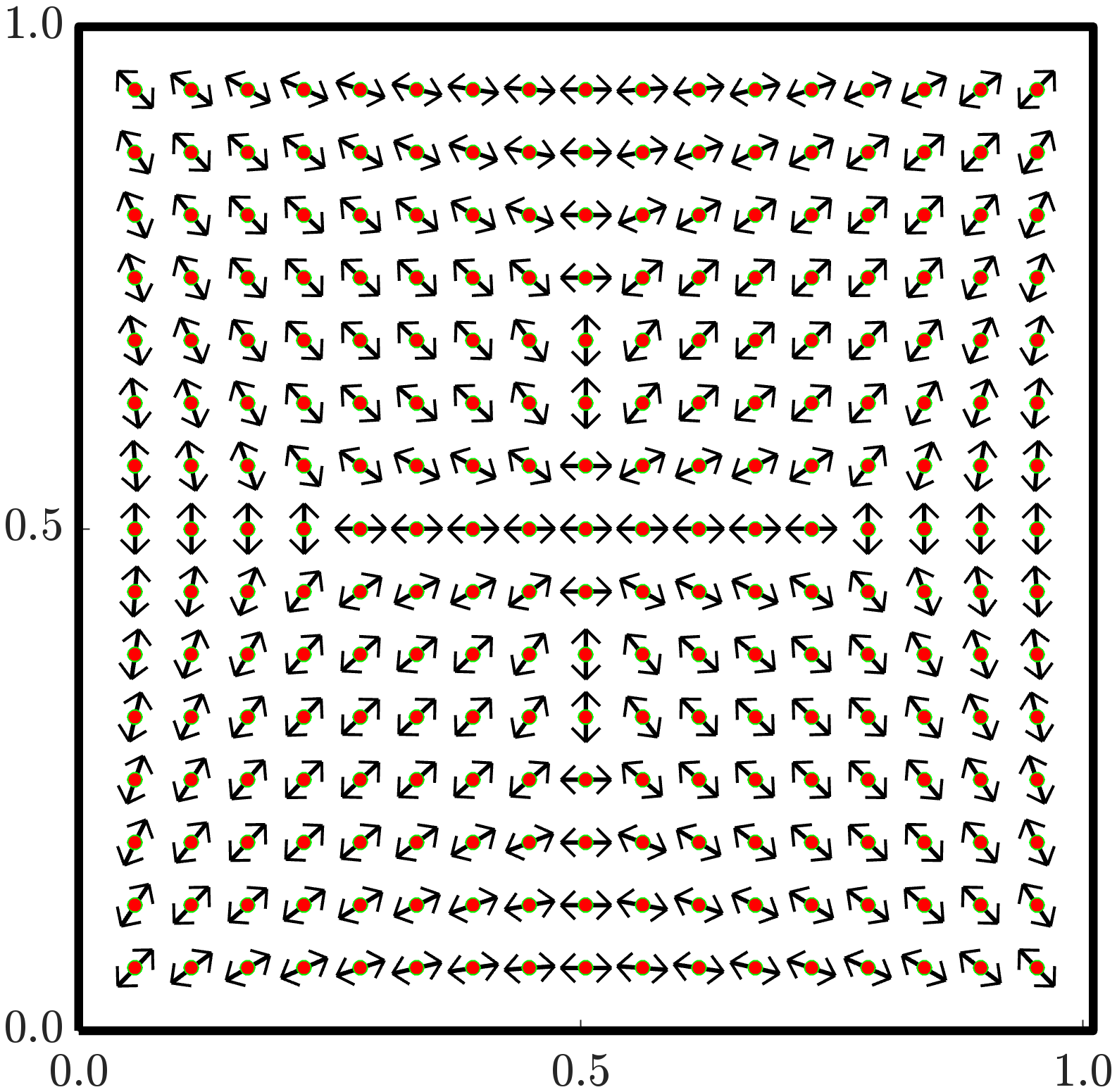}}
\hspace{0.5cm}
\subfigure[$L = 1.04$, $d = 1$.]{\includegraphics[height = 1.8in]{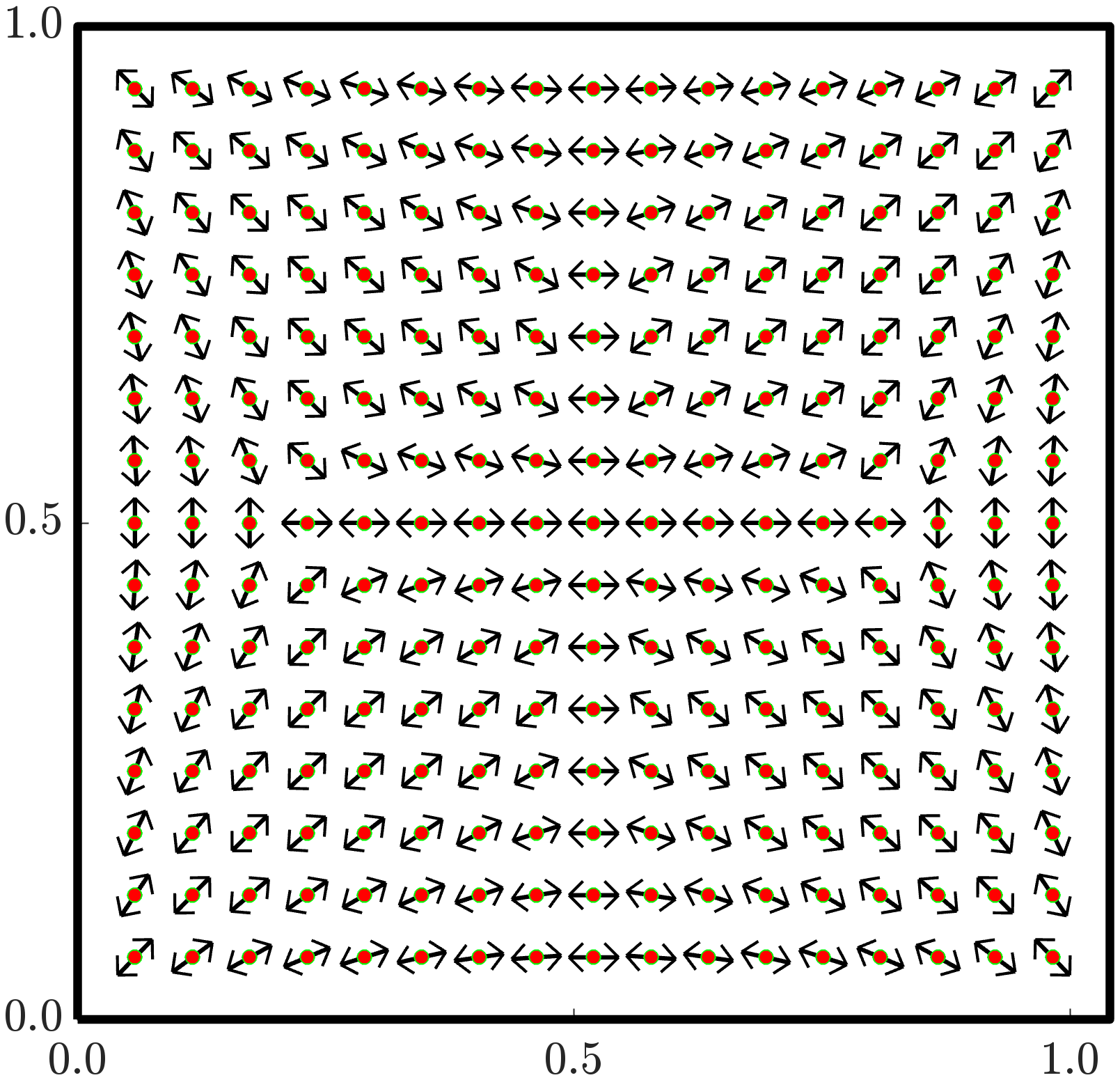}}\\
\subfigure[$L = 1.1$, $d = 1$.]{\includegraphics[height = 1.8in]{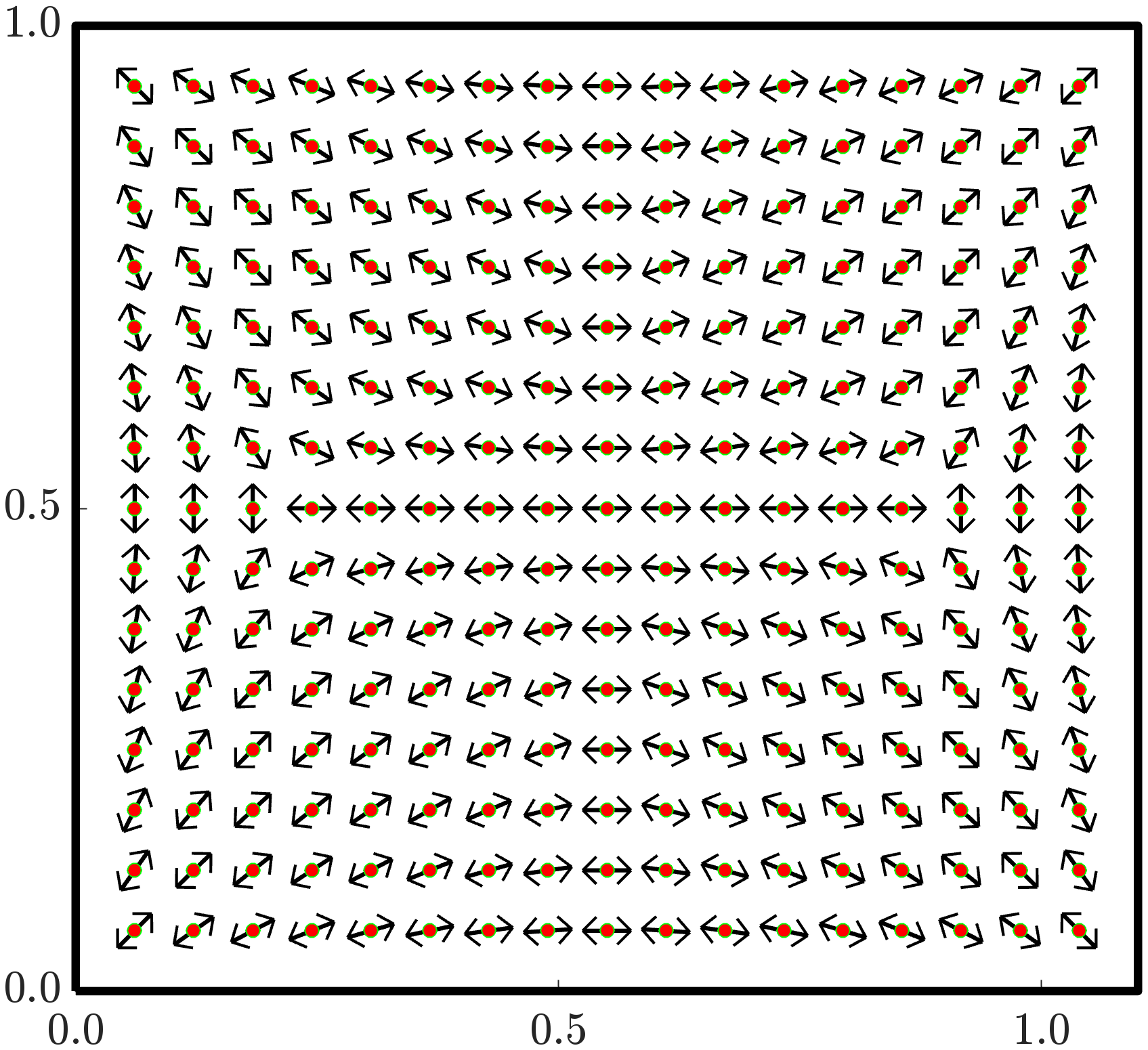}}\hspace{1cm}
\subfigure[$L = 1.5$, $d = 1$.]{\includegraphics[height = 1.8in]{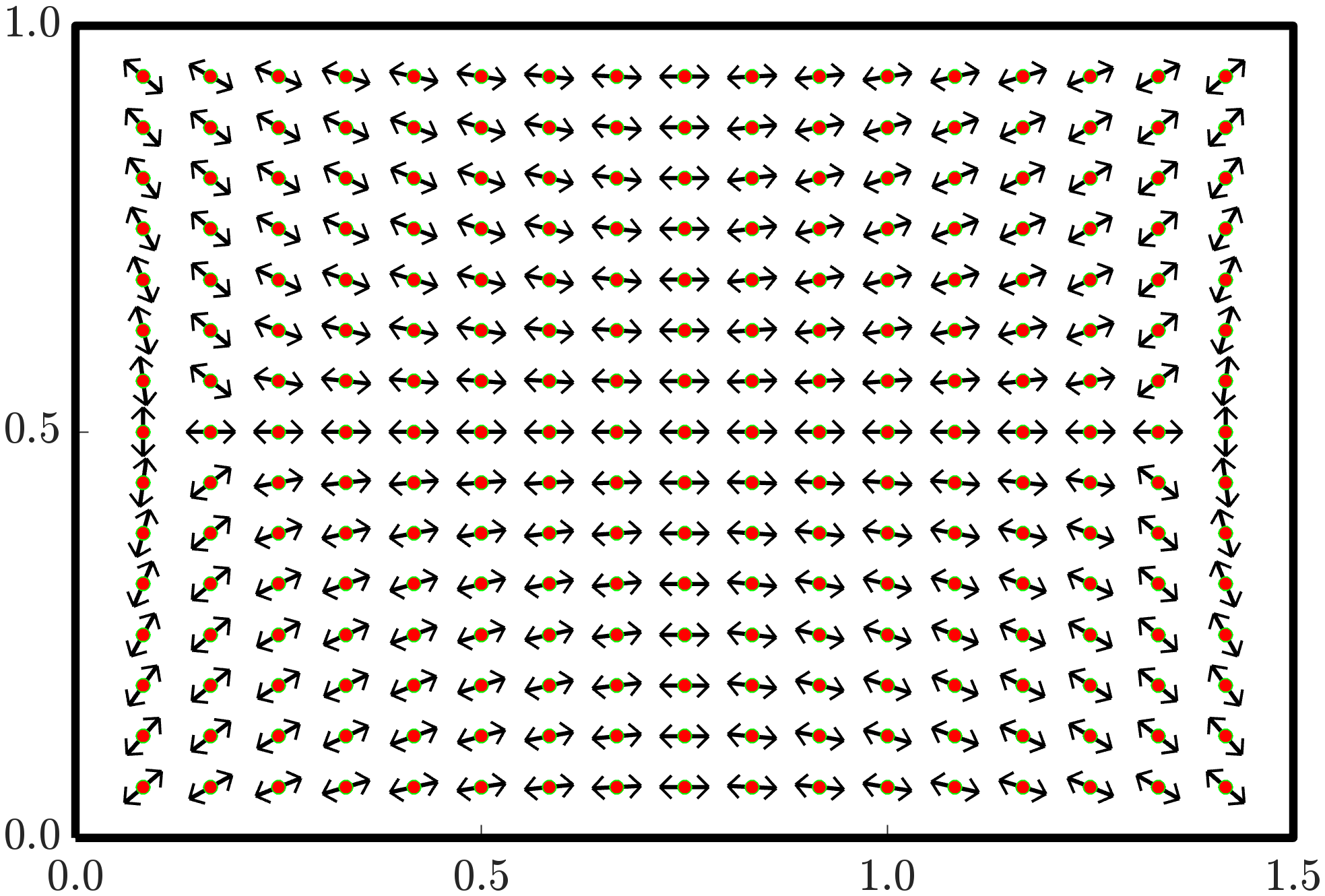}}
\caption{Minimization of $\tau_2$ for a single elliptical trap placed in the rectangular domain $\Omega = [0,L]\times[0,d]$ for $d = 1$ and $L = 1$ (a), $L=1.01$ (b) $L=1.02$ (c) $L=1.04$ (d) $L=1.1$ (e) $L=1.5$. The directional arrow indicates the direction on which the semi-major axis should be aligned to minimize $\tau_2$, the higher order GMFPT correction term.\label{fig:rectangle_opt}}
\end{figure}

\begin{figure}[htbp]
\centering
\subfigure[$a=1,\ b=1$.]{\includegraphics[height = 1.6in]{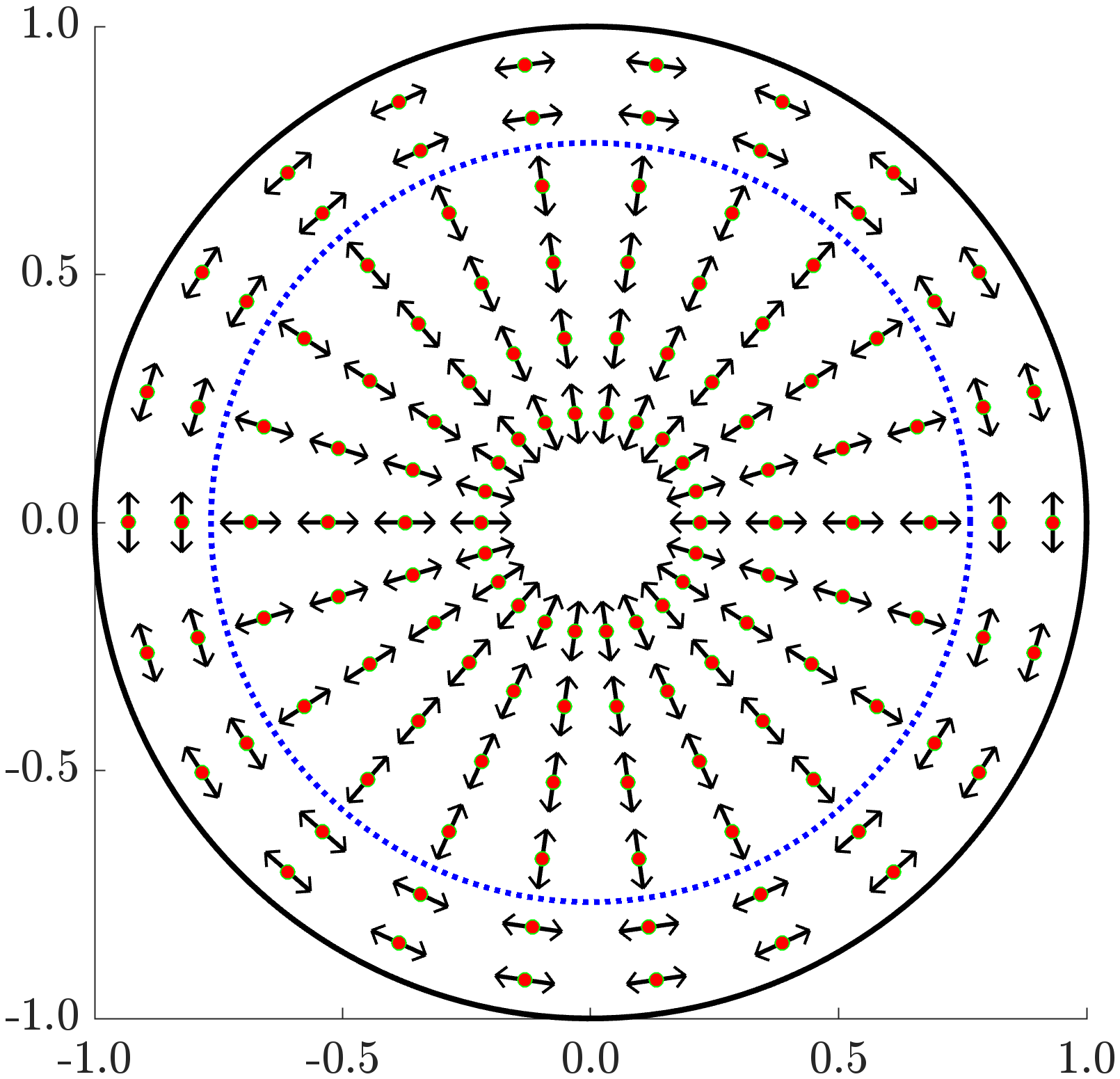}}\hspace{0.5cm}
\subfigure[$a=1.1,\ b=1$.]{\includegraphics[height = 1.6in]{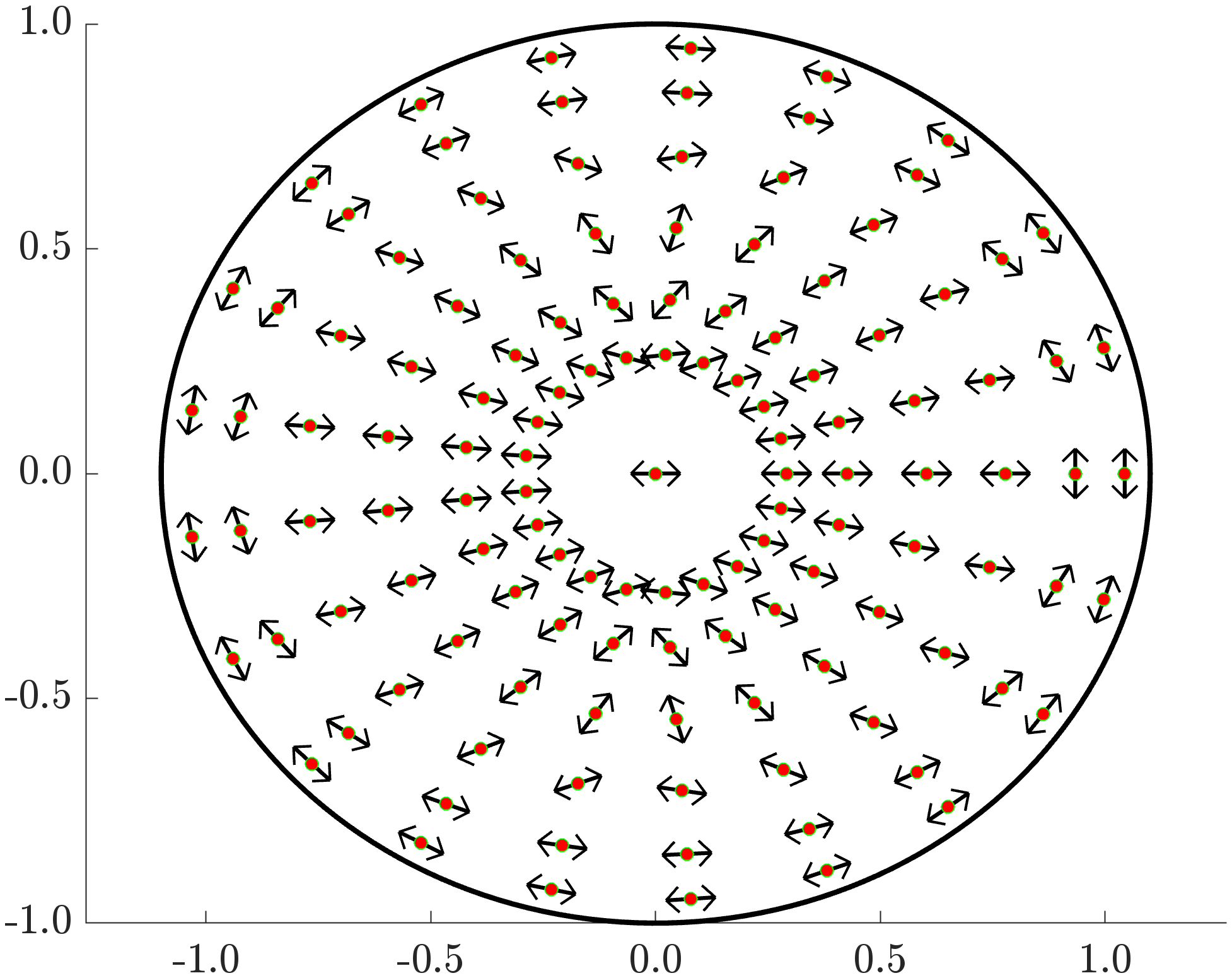}}\hspace{0.5cm}
\subfigure[$a=1.5,\ b=1$.]{\includegraphics[height = 1.6in]{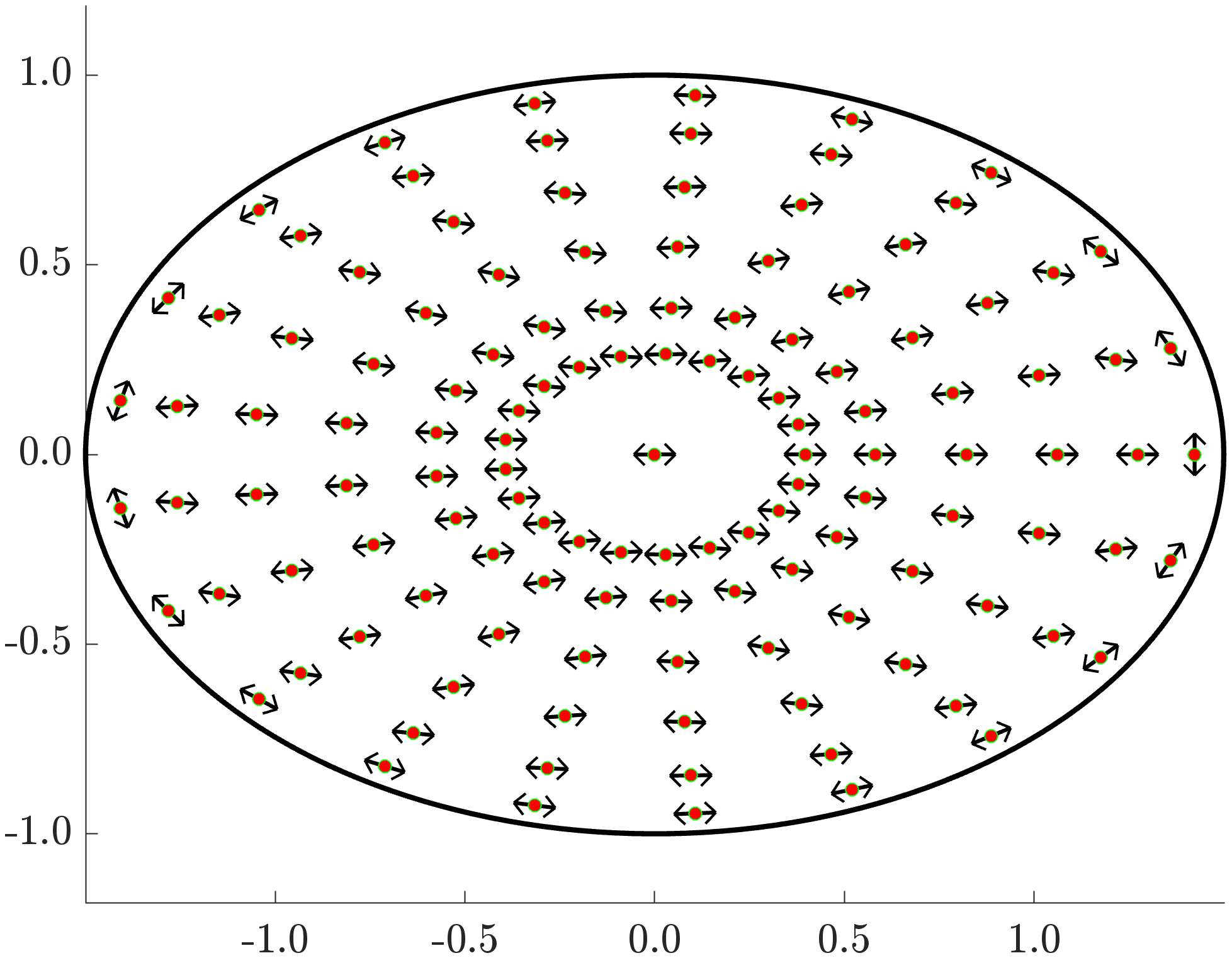}}
\caption{Minimization of $\tau_2$ for circular (a) and elliptical domains (b-c) at various locations. The directional arrow indicates the direction along which the semi-major axis should be aligned so that the correction term to the GMFPT is minimized. In Panel (a), the dashed blue line is the disk of radius $r_c = \sqrt{2-\sqrt{2}} \approx 0.7654$ where the optimal orientation flips. \label{fig:ellipse_opt}}
\end{figure}

In Fig.~\ref{fig:rectangle_opt} and Fig.~\ref{fig:ellipse_opt}, we plot vector fields of the orientation vector $\textbf{p}$ for the example of rectangular and elliptical domains, respectively. In each case, we form the necessary derivatives of $R(\bx;\bxi)$ by combining the rapidly convergent series stated in Appendix \ref{sec:Greens_functions} together with centered finite difference approximations for the first and second derivatives. As the disk domain deforms into an ellipse, the bifurcation of the minimizing orientation noted in Fig.~\ref{fig:optim_g} is smoothed out. However, a generic observation remains that for traps centered close to a smooth boundary, $\tau_2$ is minimized by orienting the semi-major axis of the trap parallel to $\partial\Omega$. 

In the case of the rectangular domain $[0,L]\times[0,d]$, we observe similar discontinuous structures that deform from the square case ($L=d$) as the rectangle elongates $(L>d)$. Curiously, the minimizing orientation of $\tau_2$ at the corners is observed to be when the semi-major axis is aligned into the corners.

\begin{figure}[htbp]
    \centering
\includegraphics[width=0.95\textwidth]{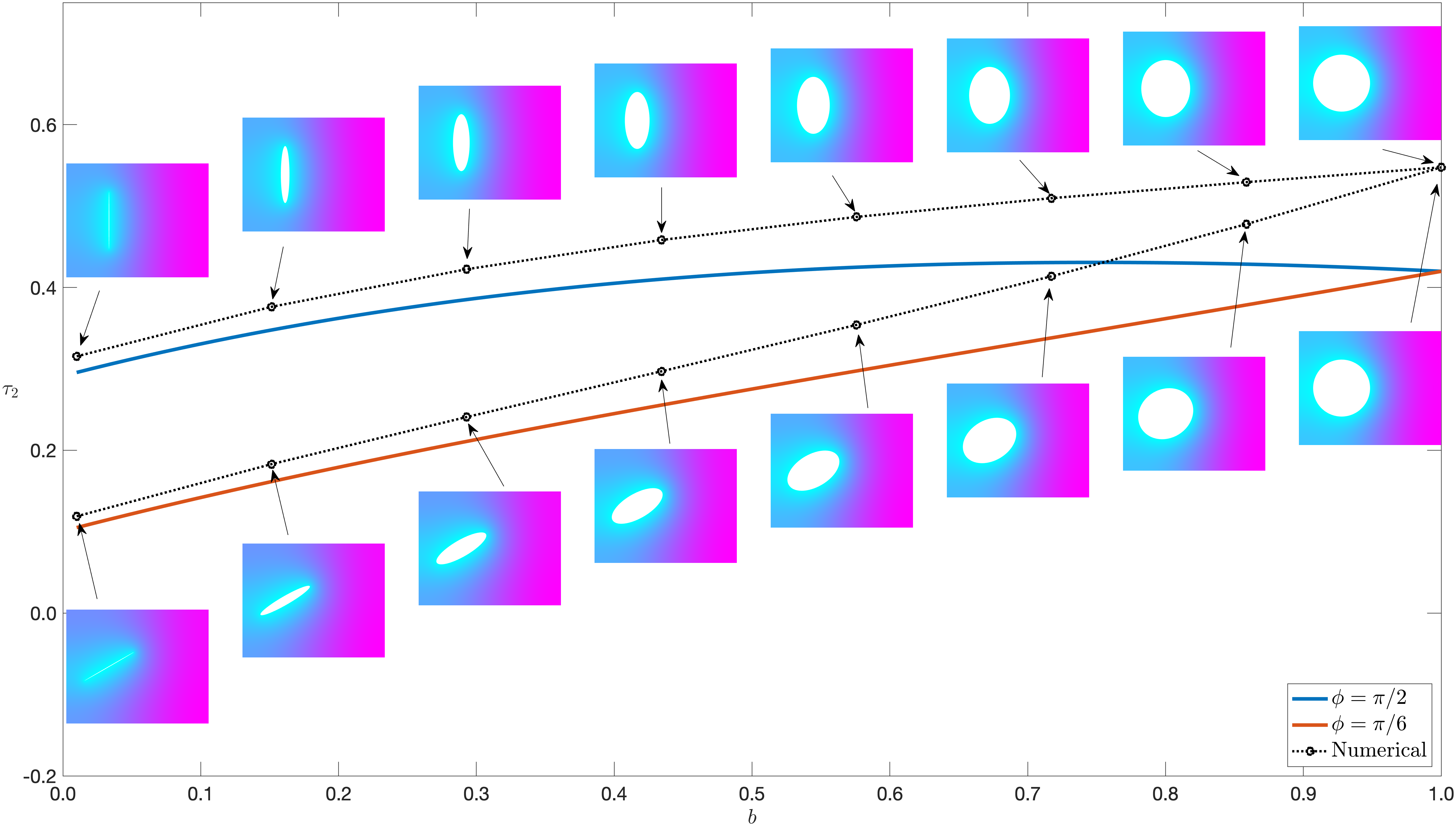}    
\caption{The effects of trap orientation and ellipticity on the high order correction to the GMFPT in the limit as $b\to0$. The correction $\tau_2 = \eps^{-2}(\tau-\tau_0)$ to the GMFPT for a rectangular domain $\Omega = [0,1]\times[0,0.8]$ with a single trap of extent $\eps = 0.2$, semi-major axis $a=1$ centered at $\bxi=[0.3,0.4]$ and varying 
semi-minor axis $b$. Curves shown for orientations $\phi =\pi/2$ and $\phi = \pi/6$ which coincide for circular traps ($b=1$).\label{fig:ellipseSlit}}
\end{figure}

\begin{exmp}
The limit of an infinitely thin ellipse to a slit.
\end{exmp} 
In this example we consider the GMFPT in the limit as $b\to0$ as the elliptical trap tends towards a thin slit. The formula for the GMFPT is uniformly valid in this limit and we find from \eqref{eqn:GlobalTauSimp} that
\begin{equation}
    \lim_{b\to0^+}\tau = \tau_0 + \eps^2 \left[ \frac{a^2}{4} -  \frac{\pi a^2|\Omega|}{2} (R_{\xi_1}^2 + R_{\xi_2}^2) + \frac{a^2|\Omega| }{4} \bp \cdot \begin{bmatrix}\cos 2\phi \\ \sin 2\phi \end{bmatrix} \right].
\end{equation}
In the case of a rectangular domain, we use this formula to explore the effect of trap orientation on the GMFPT. In Fig.~\ref{fig:ellipseSlit} we display $\tau_2$ for a trap with extent $\eps=0.2$, centered at $\bxi = [0.3,0.4]$ inside the rectangular domain $\Omega = [0,1]\times[0,0.8]$. The two curves plotted are for trap orientations $\phi = \{\pi/2,\pi/6\}$ and for $a=1$ and varying $0<b<1$. As expected, we see no effects of orientation when $a=b=1$ and smooth behavior as $b\to0$.

\section{Discussion}\label{sec:discusson}

In this work we have developed a high order matched asymptotic expansion for the MFPT of a Brownian walker to a small trap, enclosed in a two dimensional domain $\Omega$. The high order correction term describes the effect that the orientation of the trap has on the capture rate. We investigated the role that trap orientation has on the MFPT and the GMFPT, observing a sensitive dependence on the centering point of the trap in the domain. In the specific case where the enclosing domain is a disk, we found a bifurcation where the correction to the GMFPT is minimized by orientating the semi-major axis of the trap in the radial direction when the centering point satisfies $0<|\bxi|<\sqrt{2-\sqrt{2}}$, and is minimized by orientating in the angular direction when $\sqrt{2-\sqrt{2}}<|\bxi|<1$ (Fig.~\ref{fig:optim_g}). The discontinuous nature of this transition appears to be related to the symmetries of the domain and similar effects are noted in rectangular domains (Fig.~\ref{fig:rectangle_opt}). In domains with smooth boundaries, such as ellipses, we observe that the discontinuity in the vector field of optimal directions is smoothed out (Fig.~\ref{fig:ellipse_opt}). However, we observe that the GMFPT correction is generally minimized when the semi-major axis of the trap is aligned parallel to the boundary.

Our exposition of the role of trap orientation has focused on the case of an elliptical trap. This simple geometry allowed for explicit calculation of several key quantities, in particular the logarithmic capacitance $d_c$, the quadrupole matrix $\mathcal{Q}$ and the moment polarization tensor $\mathcal{M}$. The determination of these key quantities was facilitated by a complex variable approach that utilized the known mapping between the disk and ellipses to solve three variations of Laplace problems (Appendix \ref{sec:disk_problems}). However, our Principal Result \eqref{eq:introExpansion}, is valid for any trap geometry that features two lines of symmetry, for example those with rectangular or dumbbell shape. To apply the present results to these geometries, it would be necessary to calculate the three previously mentioned key quantities, either by a suitable complex transformation \cite{KTW2005} or numerical method \cite{Baddoo2021}. We hypothesize that the correction term derived in \eqref{eq:introExpansion} would vanish for trap geometries with additional lines of symmetries, such as equilateral triangles or square. In such cases, further corrections would be necessary to describe the effects of orientation.

In future work, we aim to derive the MFPT in the presence of multiple absorbing elliptical bodies. This situation is relevant in the consideration of fly muscle cells which contain multiple well spaced nuclei of elliptical shape and correlated orientations \cite{WINDNER2019}. This work shows the steps required to derive higher order corrections in a variety of narrow capture problems.

\section*{Acknowledgments} AEL acknowledges support from the NSF under award DMS 2052636.

\appendix

\section{Neumann Green's functions for disks, rectangles and ellipses}\label{sec:Greens_functions}
Here we state some known expressions for the Neumann Green's function \eqref{neum_g} for the disk, ellipse and rectangle domains.

\subsection{Neumann Green's function for a disk}\label{sec:greens_appendix}
In the case of the disk domain $\Omega = \{ \bx = (x_1,x_2) \ | \ x_1^2 + x_2^2 \leq1\}$ and source $\bxi = (\xi_1,\xi_2)$, we have \cite{KTW2005} that
\[
G = \frac{-1}{2\pi} \log|\bx-\bxi| + R(\bx;\bxi); \qquad R(\bx;\bxi) = -\frac{1}{2\pi} \left[ \frac12 \log(1+ |\bx|^2 |\bxi|^2 - 2 \bx\cdot \bxi) - \frac{1}{2}(|\bxi|^2 + |\bx|^2) + \frac{3}{4} \right].
\]
We calculate gradients $\nabla_{\bx} = (\partial_{x_1},\partial_{x_2})$, $\nabla_{\bxi} = (\partial_{\xi_1},\partial_{\xi_2})$ as
\bsub\label{eq:AppGreensDisk}
\begin{gather}
\label{eq:AppGreensDisk_a} \nabla_{\bx}G = \frac{-1}{2\pi} \frac{\bx - \bxi}{|\bx-\bxi|^2}  + \nabla_{\bx}R ; \qquad \nabla_{\bx}R = \frac{-1}{2\pi}\left[ \frac{|\bxi|^2 \bx - \bxi}{1+ |\bx|^2 |\bxi|^2 - 2 \bx\cdot \bxi} - \bx \right],\\
 \nabla_{\bxi}G = \frac{1}{2\pi} \frac{\bx - \bxi}{|\bx-\bxi|^2}  +  \nabla_{\bxi}R; \qquad \nabla_{\bxi}R = \frac{-1}{2\pi}\left[ \frac{|\bx|^2 \bxi - \bx}{1+ |\bx|^2 |\bxi|^2 - 2 \bx\cdot \bxi} - \bxi \right].
\end{gather}
As $\bx\to\bxi$ we have
\[
R(\bxi;\bxi) =\frac{-1}{2\pi} \left[ \log(1- |\bxi|^2) - |\bxi|^2 + \frac34 \right], \qquad \nabla_{\bxi}R(\bxi;\bxi) = \frac{1}{2\pi}\left[ \frac{2-|\bxi|^2}{1-|\bxi|^2} \right] \bxi .
\]
The second derivatives for the Hessian are
\begin{gather}
\label{eq:AppGreensDisk_b} \frac{\partial^2 R}{\partial{\xi_j}^2} = \frac{-1}{2\pi} \left[ \frac{|\bx|^2(1+ |\bx|^2 |\bxi|^2 - 2\bx\cdot \bxi) - 2(|\bx|^2 \xi_j - x_j)^2  }{(1+ |\bx|^2 |\bxi|^2 - 2 \bx\cdot \bxi)^2} -1\right],\quad \frac{\partial^2 R}{\partial{\xi_1\partial \xi_2}} = \frac{1}{\pi} \frac{(|\bx|^2 \xi_1-x_1)(|\bx|^2 \xi_2 - x_2)}{(1+ |\bx|^2 |\bxi|^2 - 2 \bx\cdot \bxi)^2},\\[4pt]
\label{eq:AppGreensDisk_c} \frac{\partial^2 R}{\partial{x_j}^2} = \frac{-1}{2\pi} \left[ \frac{|\bxi|^2(1+ |\bx|^2 |\bxi|^2 - 2\bx\cdot \bxi) - 2(|\bxi|^2 x_j - \xi_j)^2  }{(1+ |\bx|^2 |\bxi|^2 - 2 \bx\cdot \bxi)^2} -1\right],\quad \frac{\partial^2 R}{\partial x_1 \partial x_2} = \frac{1}{\pi} \frac{(|\bxi|^2 x_1-\xi_1)(|\bxi|^2 x_2 - \xi_2)}{(1+ |\bx|^2 |\bxi|^2 - 2 \bx\cdot \bxi)^2}.
\end{gather} 
The terms $(R_{\xi_1\xi_1} - R_{\xi_2 \xi_2})$ and $R_{\xi_1\xi_2}$ as $\bx\to\bxi$ are then
\begin{equation}
\lim_{\bx\to\bxi} \left( \frac{\partial^2 R}{\partial{\xi_1}^2}-\frac{\partial^2 R}{\partial{\xi_2}^2} \right) = \frac{1}{\pi} \frac{\xi_1^2 - \xi_2^2}{(1-|\bxi|^2)^2}, \qquad \lim_{\bx\to\bxi} \frac{\partial^2 R}{\partial \xi_1 \partial \xi_2} = \frac{1}{\pi} \frac{\xi_1\xi_2}{(1-|\bxi|^2)^2}.
\end{equation}

\esub

\subsection{Neumann Green's function for a rectangle}\label{sec:GreensRect}

The Green's function for a rectangle $\Omega = [0,L]\times[0,d]$ is known \cite{WWK09,ChenWard2011} in the form of a rapidly convergent series. For $\bx = (x_1,x_2)$ and $\bxi = (\xi_1,\xi_2)$ we have
\bsub
\begin{align}
\nonumber R(\bx;\bxi) &= \frac{-1}{2\pi} \sum_{n=0}^{\infty} \log( |1-q^n z_{+,+}||1-q^n z_{+,-}||1-q^n z_{-,+}||1-q^n \zeta_{+,+}||1-q^n \zeta_{+,-}||1-q^n \zeta_{-,+}||1-q^n \zeta_{-,-}|)\\
{}&  -\frac{1}{2\pi} \log \frac{|1- z_{-,-}|}{|r_{-,-}|} + \frac{L}{d}\left[ \frac13 - \frac{\max(x_1,\xi_1)}{L} + \frac{x_1^2+\xi_1^2}{2L^2}\right] - \frac{1}{2\pi} \sum_{n=1}^{\infty} \log| 1- q^n z_{-,-}|.
\end{align}
where 
\begin{align}
z_{\pm,\pm} &\equiv e^{\mu r_{\pm,\pm}/2}, \qquad \zeta_{\pm,\pm} \equiv e^{\mu \rho_{\pm,\pm}/2}, \qquad \mu \equiv \frac{2\pi}{d}, \quad q \equiv e^{-\mu L},\\
r_{+,\pm} &\equiv -|x_1 + \xi_1| + i (x_2 \pm \xi_2), \qquad r_{-,\pm} \equiv -|x_1 - \xi_1| + i (x_2 \pm \xi_2),\\[5pt]
\rho_{+,\pm} &\equiv |x_1 + \xi_1| + i (x_2 \pm \xi_2) -2L, \qquad \rho_{-,\pm} \equiv |x_1 - \xi_1| + i (x_2 \pm \xi_2) - 2L.
\end{align}
\esub
The self interaction term $R(\bx;\bx)$ is given by 
\bsub
\begin{align}
\nonumber R(\bx;\bx) &= \frac{-1}{2\pi} \sum_{n=0}^{\infty} \log( |1-q^n z^0_{+,+}||1-q^n z^0_{+,-}||1-q^n z^0_{-,+}||1-q^n \zeta^0_{+,+}||1-q^n \zeta^0_{+,-}||1-q^n \zeta^0_{-,+}||1-q^n \zeta^0_{-,-}|)\\
{}& + \frac{L}{d} \left( \frac13 - \frac{x_1}{L} + \frac{x_1^2}{L^2} \right) -\frac{1}{2\pi}\log \left( \frac{\pi}{d}\right)  - \frac{1}{2\pi}\sum_{n=1}^{\infty} \log(1- q^n ).
\end{align}
where
\begin{align}
z^0_{+,+} &\equiv e^{\mu(-x_1 + i x_2)}, \qquad z^0_{+,-} \equiv e^{-\mu x_1 }, \qquad z^0_{-,+} \equiv e^{\mu i x_2 }, \\[5pt]
\zeta^0_{+,+} & \equiv e^{\mu (x_1 -L + i x_2)}, \qquad \zeta^0_{-,+} \equiv e^{\mu (-L + i x_2)}, \qquad \zeta^0_{+,-} \equiv e^{\mu (x_1 -L)}, \qquad \zeta^0_{-,-} \equiv e^{- \mu L}.
\end{align}
\esub

\subsection{Neumann Green's function for an ellipse}\label{sec:GreensEllipse}

A rapidly convergent series for the solution of \eqref{neum_g} in the elliptical domain $\Omega = \{ \bx = (x_1,x_2) \; | \; (x_1/a)^2 + (x_2/b)^2 \leq 1\}$ was derived in \cite{Iyaniwura2021_b}. For completeness we restate the final result here. The first step is to introduce the transformation,
\bsub\label{eq:transform}
\begin{equation}
x_1 = f \cosh \xi \cos \eta, \qquad x_2 = f \sinh \xi \sin \eta, \qquad f = \sqrt{a^2-b^2},
\end{equation}
which maps $\bx = (x_1,x_2)\in\Omega$ to the rectangle $0\leq \xi \leq \xi_b$ and $0\leq \eta \leq2\pi$ where $a = f \cosh\xi_b$ and $b = f\sinh\xi_b$ so that
\begin{equation}
f = \sqrt{a^2-b^2}, \qquad \xi_b = \tanh^{-1} \frac{b}{a} = -\frac{1}{2} \log \gamma, \qquad \gamma = \left(\frac{a-b}{a+b}\right).
\end{equation}
For a pair $(x_1,x_2)$, the corresponding $(\xi,\eta)$ satisfy
\[
\xi = \frac12 \log\Big(1-2s + s\sqrt{s^2-s} \Big), \qquad s = \frac{-\mu - \sqrt{\mu^2 + 4 f^2 y^2}}{2f^2}, \qquad \mu = x_1^2 + x_2^2 - f^2.
\]
For $\eta_{\ast} = \sin^{-1} (\sqrt{p})$, the value of $\eta$ is given by
\begin{equation}
    \eta = \left\{ \begin{array}{lr}
\eta_{*}, & \text{for } x_1 \geq0,\ x_2 \geq0\\
\pi-\eta_{*}, &  \text{for } x_1 <0,\ x_2 \geq0\\
\pi+\eta_{*}, &  \text{for } x_1 \leq0,\ x_2 <0\\
2\pi-\eta_{*}, &  \text{for } x_1 >0,\ x_2 <0\\        
    \end{array} \right. \, , \qquad \text{where}
    \qquad p  = \frac{-\mu + \sqrt{\mu^2 + 4f^2 y^2}}{2f^2}.
\end{equation}
\esub
For points $\bx= (x_1,x_2)$ and $\by = (y_1,y_2)$, the Green's function $G(\bx;\by)$ for $\bx\neq\by$ is given by
\begin{equation}
    G(\bx;\by) = \frac{1}{4|\Omega|}(|\bx|^2 + |\by|^2) - \frac{3}{16|\Omega|}(a^2 + b^2) - \frac{1}{4\pi}\log\gamma - \frac{1}{2\pi}\max(\xi,\xi_0) - \frac{1}{2\pi}\sum_{n=0}^{\infty} \log \left( \prod_{j=1}^8 |1- \gamma^{2n} z_j| \right),
\end{equation}
where $|\Omega| = \pi ab$. The complex constants $z_1,\ldots, z_8$ are defined in terms of $(\xi,\eta)$, $(\xi_0,\eta_0)$ and $\xi_b$ by
\begin{gather*}
    z_1 = e^{-|\xi-\xi_0| + i (\eta-\eta_0)}, \quad z_2 = e^{|\xi-\xi_0| - 4\xi_b + i (\eta-\eta_0)}, \quad z_3 = e^{(\xi+\xi_0) - 2\xi_b + i (\eta-\eta_0)},\\[4pt]
    z_4 = e^{(\xi+\xi_0) - 2\xi_b + i (\eta-\eta_0)}, \quad z_5 = e^{(\xi+\xi_0) - 4\xi_b + i (\eta+\eta_0)}, \quad z_6 = e^{-(\xi+\xi_0) + i (\eta+\eta_0)},\\[4pt]
    z_7 = e^{|\xi+\xi_0| - 2\xi_b + i (\eta+\eta_0)}, \quad z_8 = e^{-|\xi+\xi_0| - 2\xi_b + i (\eta+\eta_0)}.
\end{gather*}
The point $(x_1,x_2)$ is mapped to $(\xi,\eta)$ while the source point $(y_1,y_2)$ is mapped to $(\xi_0,\eta_0)$ by the transformation \eqref{eq:transform}. The quantity $R(\by;\by)$ is given by
\begin{align*}
    \nonumber R(\by;\by) &= \frac{|\by|^2}{2|\Omega|} - \frac{3}{16|\Omega|}(a^2 + b^2) + \frac{1}{2\pi} \log(a+b) - \frac{\xi_0}{2\pi} + \frac{1}{4\pi} \log\big( \cosh^2 \xi_0 - \cos^2 \eta_0 \big)\\[4pt]
    &- \frac{1}{2\pi} \sum_{n=1}^{\infty} 
    \log (1- \gamma^{2n}) - \frac{1}{2\pi} \sum_{n=0}^{\infty} \log \left( \prod_{j=2}^8 |1-\gamma^{2n} z_j^0| \right)
\end{align*}
Here the constants $z_j^0$ for $j=2,\ldots,8$ are 
\begin{gather*}
    z_2^0 = \gamma^2, \quad z_3^0 = \gamma e^{-2\xi_0}, \quad z_4^0 = \gamma e^{2\xi_0}, \quad z_4^0 = \gamma^2 e^{2\xi_0+2i\eta_0},\\[4pt]
    z_6^0 = e^{-2\xi_0 + 2i \eta_0}, \quad z_7^0 = \gamma e^{2i\eta_0}, \quad z_8^0 = \gamma e^{2i \eta_0}, \qquad \gamma = \frac{a-b}{a+b}.
\end{gather*}

\section{Inner problems for the exterior of the ellipse}\label{sec:disk_problems}

We solve a variety of Laplace equations posed in the exterior of the elliptical domain $\A = \{ \by = (y_1,y_2) \ | \ y_1^2/a^2 + y_2^2/b^2 < 1  \}$. In each case, we make use of the complex transformation
\begin{equation}\label{app:eq1}
\by = \alpha \, \bz + \frac{\beta}{\bz},
\end{equation}
which maps the unit disk to the ellipse $\A$ with semi-major and semi-minor axes $a$ and $b$ respectively with the semi-major axis aligned on the horizontal axis. On the unit disk $\bz = e^{i\theta}$, we have that
\[
\by = \alpha  e^{i\theta} + \beta  e^{-i\theta} = (\alpha + \beta) \cos\theta + i(\alpha - \beta)\sin\theta = a \cos\theta + i b \sin\theta.
\]
The mapping parameters $\alpha$ and $\beta$ are then
\begin{equation}\label{app:alfbet}
\alpha = \frac{a+b}{2}, \qquad \beta = \frac{a-b}{2}.
\end{equation}
The Laplace equations to be considered will be solved in the unit disk then mapped to the ellipse by the inverse transformation of \eqref{app:eq1}. The large argument behavior of the inverse transform of \eqref{app:eq1} is calculated as follows for $|\by|\gg1$
\begin{equation}\label{app:eq2} 
 \bz \sim \frac{\by}{\alpha} - \frac{\beta}{\by} =  \frac{\by}{\alpha} - \frac{\beta\bar{\by}}{|\by|^2}, \qquad 
 |\bz| = \frac{|\by|}{\alpha}\left[1 - \frac{\alpha\beta}{|\by|^4} \by^T \begin{pmatrix} 1 & \phantom{-}0\\ 0& -1\end{pmatrix} \by + \bigoh(|\by|^{-4}) \right] , \qquad |\by| \to\infty.
\end{equation}

\subsection{The order $\bigoh(\eps^0)$ problem}\label{sec:Quadterm}

The leading order inner problem is given by
\bsub\label{vc}
\begin{gather}
\label{vc_a} \Delta_\by v_{0c}  = 0\,, \qquad  \by\in\mathbb{R}^2\setminus \A; \qquad v_{0c}  = 0, \qquad \by \in \partial\A;
\\[5pt]
\label{vc_b}   v_{0c} = \log |\by| - \log d_c  + \frac{\mathbf{d} \cdot \by}{|\by|^2} + \frac{\by^{T} \tilde{\Q} \by}{|\by|^4} +  \cdots \, \quad   |\by|\to\infty.
\end{gather}
\esub
In classic potential theory, the logarithmic term is the monopole, $\mathbf{d}$ is the dipole vector and $\tilde{\Q}$ is the quadrupole matrix.
Our goal is to obtain the solution of \eqref{vc} and identify the logarithmic capacitance $d_c$ and quadrupole matrix $\tilde{\Q}$. In the scenario of the unit disk ($a=b=1$), the unique solution is $\log|\bz|$. For the problem \eqref{vc}, we apply the mapping \eqref{app:eq2} to obtain the far field behavior
\begin{equation}\label{app:eq3}
v_{0c}(\by) = \log|\bz| \sim \log|\by| - \log \alpha - \frac{\alpha\beta}{|\by|^4} \by^T \begin{pmatrix} 1 & \phantom{-}0\\ 0& -1\end{pmatrix} \by + \bigoh(|\by|^{-4}).
\end{equation}
Hence in comparing \eqref{vc_b} with \eqref{app:eq3}, and applying \eqref{app:alfbet}, we establish that
\begin{equation}\label{app:eq4}
d_c = \log\frac{a+b}{2}, \qquad \mathbf{d} = \mathbf{0}, \qquad  \tilde{\Q}  = -\alpha\beta \begin{pmatrix} 1 & \phantom{-}0\\ 0& -1\end{pmatrix}= - \frac{a^2 -b^2}{4} \begin{pmatrix} 1 & \phantom{-}0\\ 0& -1\end{pmatrix}.
\end{equation}
The vanishing dipole vector $\mathbf{d} = 0$ is consistent with the two lines of symmetry of the elliptical domain. Extending this analogy, domains with four lines of symmetries (e.g. square) would have a vanishing quadrupole matrix requiring an even higher order of expansion.

\subsection{The order $\bigoh(\eps^1)$ problem}\label{app:orderone}
The first order inner problem is vector valued and given by
\bsub\label{vc1}
\begin{gather}
\Delta_\by \textbf{v}_{1c} = 0\,, \quad  \by\in\mathbb{R}^2\setminus \A; \qquad  \textbf{v}_{1c}  = 0, \quad \by \in \partial\A;
\\[5pt]
 \textbf{v}_{1c} = \by + \frac{ \tilde{\M} \by}{|\by|^2} + \cdots \, \quad   |\by|\to\infty.
\end{gather}
\esub
For the case of the unit disk, we have that $\textbf{v}_{1d}= \bz - \bz/|\bz|^2$. Hence for the elliptical domain under the transformation \eqref{app:eq2}, we have that
\begin{equation}
\textbf{v}_{1c}(\by) = \alpha \textbf{v}_{1d}(\by) = \by + \frac{\tilde{\M}\by}{|\by|^2} +\cdots \qquad |\by|\to\infty; \qquad \tilde{\M} = - \alpha\begin{pmatrix} a & 0\\ 0& b\end{pmatrix}.
\end{equation}
\subsection{The order $\bigoh(\eps^2)$ problem}\label{app:eqnv2}
In this subsection, we derive the solution of the higher order correction problem
\bsub\label{vc2}
\begin{gather}
\label{vc2_a} \Delta_\by v_{2c} = -1, \qquad \by \in \mathbb{R}^2\setminus\A; \qquad  v_{2c}  = 0,  \quad \by\in \partial\A;\\[4pt]
\label{vc2_b} v_{2c} = -\frac{|\by|^2}{4} + \by^{T} \mathcal{B} \by + \cdots \, \quad   |\by|\to\infty, \qquad \mathcal{B} = \begin{bmatrix} \mathcal{B}_{11} &\phantom{-}\mathcal{B}_{12} \\ \mathcal{B}_{12} & -\mathcal{B}_{11}\end{bmatrix}.
\end{gather}
\esub
where $\mbox{Trace}(\mathcal{B}) = 0$. The general solution takes the form $v_{2c} = -\frac{1}{4}|\by|^2 + v_{2h}$ where the homogeneous solution satisfies
\bsub\label{vh2}
\begin{gather}
\label{vh2_a} \Delta_{\by} v_{2h} = 0, \qquad \by \in \mathbb{R}^2\setminus\A; \qquad v_{2h} = \frac{|\by|^2}{4}, \qquad \by\in \partial\A;\\[4pt]
\label{vh2_b} v_{2h} \sim  \by^{T} \mathcal{B} \by + \cdots \quad |\by|\to\infty.
\end{gather}
\esub
As with previous solutions of inner problems, we solve the corresponding problem on the disk and use the complex transformation \eqref{app:eq1} to map the solution to $\A$. Recalling that when $\bz=e^{i\theta}$, we have $|\by|^2 = (\alpha^2 + \beta^2) + 2\alpha\beta \cos 2\theta$. Hence, the homogeneous solutions is expressed in terms of the disk solution in complex form as
\begin{equation}\label{eqn:v2h_gen}
v_{2h} = a_1 + \mbox{Re}[ b_1\bz^2 +b_2 \bz^{-2}] + \mbox{Im} [c_1\bz^2 +c_2 \bz^{-2}].
\end{equation}
On the boundary $\bz=e^{i\theta}$, we have the conditions
\[
 \frac14(\alpha^2 + \beta^2) + \frac12 \alpha\beta \cos 2\theta = a_1 +  (b_1 + b_2)\cos 2\theta + (c_1 + c_2)\sin 2\theta,
\]
which yields the conditions
\[
a_1 = \frac14 (\alpha^2 + \beta^2), \qquad b_1 + b_2 = \frac12\alpha\beta, \qquad c_1+c_2 = 0.
\]
To establish the behavior as $|\by| \to \infty$, we consider from \eqref{app:eq2} that
\[
\bz^2 \sim \left(\frac{\by}{\alpha} - \frac{\beta}{\by}\right)^2 = \frac{\by^2}{\alpha^2} - 2\frac{\beta}{\alpha} + \frac{\beta^2}{\by^2}, \quad \mbox{as} \quad |\by|\to\infty.
\]
The large argument behavior of $v_{2h}$ given in \eqref{eqn:v2h_gen} is
\begin{align}
    \nonumber v_{2h} &\sim \frac{b_1}{\alpha^2} \mbox{Re}[\by^2] + \frac{c_1}{\alpha^2}  \mbox{Im}[\by^2] + \frac14 (\alpha^2 + \beta^2) - 2\frac{\beta}{\alpha} b_1 + \bigoh(|\by|^{-2}).
\end{align}
Now, comparing with \eqref{vh2_b}, we see that 
\begin{equation}
b_1 = \alpha^2 \mathcal{B}_{11}, \qquad c_1 = \alpha^2 \mathcal{B}_{12}.
\end{equation}
Hence, we obtain the large argument behavior of \eqref{vh2} to be 
\bsub
\begin{align}
    v_{2h} &\sim \by^T \mathcal{B}\by + \frac{1}{4}(\alpha^2 + \beta^2) - 2\alpha\beta \, \mathcal{B}_{11} +  +  \bigoh(|\by|^{-2}),
\end{align}
\esub
and finally establish the large argument behavior of \eqref{vc2} to be
\begin{equation}\label{app:2ndorder2}
 v_{2c} = -\frac{|\by|^2}{4} + \by^{T} \mathcal{B} \by + d_{2c} + \mathcal{O}(|\by|^{-2}), \quad \mbox{as} \quad |\by|\to\infty;  \qquad d_{2c} = \frac{1}{4}( \alpha^2 + \beta^2) - 2\alpha\beta \, \mathcal{B}_{11} .
\end{equation}

\bibliographystyle{siam}
\bibliography{refs}

\end{document}